\documentclass[10pt]{amsart}
\usepackage{amsfonts}
\usepackage{amssymb}
\usepackage{amsmath}
\usepackage[dvips]{graphicx}

\input{tcilatex}

\begin{document}
\title[Certain subclasses of multivalent functions]{\textbf{Certain
subclasses of multivalent functions defined by new multiplier transformations%
}}
\author{Erhan Deniz and Halit orhan}
\address{\thinspace Department of Mathematics, Faculty of Science, Ataturk
University, 25240 Erzurum, Turkey}
\email{edeniz@atauni.edu.tr}
\email{horhan@atauni.edu.tr}
\date{}

\begin{abstract}
In the present paper the new multiplier transformations $\mathrm{{\mathcal{J}%
}}_{p}^{\delta }(\lambda ,\mu ,l)$ $(\delta ,l\geq 0,\;\lambda \geq \mu \geq
0;\;p\in \mathrm{%
\mathbb{N}
)}$ of multivalent functions is defined. Making use of the operator $\mathrm{%
{\mathcal{J}}}_{p}^{\delta }(\lambda ,\mu ,l),$ two new subclasses $\mathcal{%
P}_{\lambda ,\mu ,l}^{\delta }(A,B;\sigma ,p)$ and $\widetilde{\mathcal{P}}%
_{\lambda ,\mu ,l}^{\delta }(A,B;\sigma ,p)$\textbf{\ }of multivalent
analytic functions are introduced and investigated in the open unit disk.
Some interesting relations and characteristics such as inclusion
relationships, neighborhoods, partial sums, some applications of fractional
calculus and quasi-convolution properties of functions belonging to each of
these subclasses $\mathcal{P}_{\lambda ,\mu ,l}^{\delta }(A,B;\sigma ,p)$
and $\widetilde{\mathcal{P}}_{\lambda ,\mu ,l}^{\delta }(A,B;\sigma ,p)$ are
investigated. Relevant connections of the definitions and results presented
in this paper with those obtained in several earlier works on the subject
are also pointed out.
\end{abstract}

\keywords{Multiplier transformations; Analytic functions; Multivalent
functions; Neighborhoods and Partial sums of analytic functions; Fractional
calculus operators.}
\subjclass[2000]{Primary 30C45}
\maketitle

\noindent

\section{\ INTRODUCTION AND DEFINITIONS}

Let $\mathrm{{\mathcal{A}}}(n,p)$ denote the class of functions normalized by%
\begin{equation}
f(z)=z^{p}+\sum_{k=n+p}^{\infty }a_{k}z^{k}\text{\ \ }\left( p,n\in \mathrm{%
\mathbb{N}
}:=\{1,2,3,...\}\right)  \tag{1.1}
\end{equation}%
which are analytic and $p-valent$ in the open unit disk $\mathrm{{\mathcal{U}%
}}=\{z:\;z\in \mathrm{%
\mathbb{C}
}\;\mathrm{and\;}\left\vert z\right\vert <1\}.$

Let $f(z)$ and $g(z)$ be analytic in $\mathrm{{\mathcal{U}}}.$ Then, we say
that the function $f$ is subordinate to $g$ if there exists a Schwarz
function $w(z),$ analytic in $\mathrm{{\mathcal{U}}}$ with $w(0)=0,$ $%
\left\vert w(z)\right\vert <1$ such that $f(z)=g(w(z))\mathrm{\;(}z\in 
\mathrm{{\mathcal{U}})}.$ We denote this subordination $f\prec g\mathrm{%
\;or\;}f(z)\prec g(z)\mathrm{\;(}z\in \mathrm{{\mathcal{U}})}.$ In
particular, if the function $g$ is univalent in $\mathrm{{\mathcal{U}}}$,
the above subordination is equivalent to $f(0)=g(0),$ $f(\mathrm{{\mathcal{U}%
}})\subset g(\mathrm{{\mathcal{U}}}).$

For $f\in \mathrm{{\mathcal{A}}}(n,p)$ given by (1.1) and $g(z)$ given by%
\begin{equation}
g(z)=z^{p}+\sum_{k=n+p}^{\infty }b_{k}z^{k}\text{\ \ }\left( p,n\in \mathrm{%
\mathbb{N}
}:=\{1,2,3,...\}\right)  \tag{1.2}
\end{equation}%
their convolution (or Hadamard product), denoted by $(f\ast g),$ is defined
as%
\begin{equation}
(f\ast g)(z):=z^{p}+\sum_{k=n+p}^{\infty }a_{k}b_{k}z^{k}=:(g\ast f)(z)\text{
\ }\left( z\in \mathrm{{\mathcal{U}}}\right) .  \tag{1.3}
\end{equation}%
Note that $f\ast g\in \mathrm{{\mathcal{A}}}(n,p).$ In particular, we set%
\begin{equation*}
\mathrm{{\mathcal{A}}}(p,1):=\mathrm{{\mathcal{A}}}_{p},\;\;\mathrm{{%
\mathcal{A}}}(1,n):=\mathrm{{\mathcal{A}}}(n),\;\;\mathrm{{\mathcal{A}}}%
(1,1):=\mathrm{{\mathcal{A}}}_{1}=\mathrm{{\mathcal{A}}}.\;\;
\end{equation*}

For a function $f$ in $\mathrm{{\mathcal{A}}}(n,p),$ we define the \textit{%
multiplier transformations }$\mathrm{{\mathcal{J}}}_{p}^{\delta }(\lambda
,\mu ,l)$ as follows:

\noindent \textbf{Definition 1.1. }Let $f\in \mathrm{{\mathcal{A}}}(n,p).$
For the parameters $\delta ,\lambda ,\mu ,l\in \mathrm{%
\mathbb{R}
};$ $\lambda \geq \mu \geq 0$ and $\delta ,l\geq 0$ define the multiplier
transformations\textit{\ }$\mathrm{{\mathcal{J}}}_{p}^{\delta }(\lambda ,\mu
,l)$ on $\mathrm{{\mathcal{A}}}(n,p)$ by the following%
\begin{equation*}
\mathrm{{\mathcal{J}}}_{p}^{0}(\lambda ,\mu ,l)f(z)=f(z)
\end{equation*}%
\begin{equation*}
~(p+l)\mathrm{{\mathcal{J}}}_{p}^{1}(\lambda ,\mu ,l)f(z)=\lambda \mu
z^{2}f^{\prime \prime }(z)+\left( \lambda -\mu +(1-p)\lambda \mu \right)
zf^{\prime }(z)+\left( p(1-\lambda +\mu )+l\right) f(z)
\end{equation*}%
\begin{align}
(p+l)\mathrm{{\mathcal{J}}}_{p}^{2}(\lambda ,\mu ,l)f(z)& =\lambda \mu z^{2}[%
\mathrm{{\mathcal{J}}}_{p}^{1}(\lambda ,\mu ,l)f(z)]^{\prime \prime }+\left(
\lambda -\mu +(1-p)\lambda \mu \right) z[\mathrm{{\mathcal{J}}}%
_{p}^{1}(\lambda ,\mu ,l)f(z)]^{\prime }  \tag{1.4} \\
& +\left( p(1-\lambda +\mu )+l\right) \mathrm{{\mathcal{J}}}_{p}^{1}(\lambda
,\mu ,l)f(z)  \notag
\end{align}%
\begin{equation*}
\mathrm{{\mathcal{J}}}_{p}^{\delta _{1}}(\lambda ,\mu ,l)(\mathrm{{\mathcal{J%
}}}_{p}^{\delta _{2}}(\lambda ,\mu ,l)f(z))=\mathrm{{\mathcal{J}}}%
_{p}^{\delta _{2}}(\lambda ,\mu ,l)(\mathrm{{\mathcal{J}}}_{p}^{\delta
_{1}}(\lambda ,\mu ,l)f(z))
\end{equation*}%
for $z\in \mathrm{{\mathcal{U}}}$ and $p,n\in \mathrm{%
\mathbb{N}
}:=\{1,2,...\}.$

If $f$ is given by (1.1) then from the definition of the multiplier
transformations $\mathrm{{\mathcal{J}}}_{p}^{\delta }(\lambda ,\mu ,l),$ we
can easily see that%
\begin{equation}
\mathrm{{\mathcal{J}}}_{p}^{\delta }(\lambda ,\mu
,l)f(z)=z^{p}+\sum_{k=n+p}^{\infty }\Phi _{p}^{k}(\delta ,\lambda ,\mu
,l)a_{k}z^{k}  \tag{1.5}
\end{equation}%
where%
\begin{equation}
\Phi _{p}^{k}(\delta ,\lambda ,\mu ,l)=\left[ \frac{(k-p)(\lambda \mu
k+\lambda -\mu )+p+l}{p+l}\right] ^{\delta }.  \tag{1.6}
\end{equation}

\noindent \textbf{Remark 1.1.} It should be remarked that the operator$%
\mathrm{{\mathcal{J}}}_{p}^{\delta }(\lambda ,\mu ,l)$ is a generalization
of many other\textit{\ operators} considered earlier. In particular, for $%
f\in \mathrm{{\mathcal{A}}}(n,p)$ we have the following:

\begin{enumerate}
\item $\mathrm{{\mathcal{J}}}_{1}^{\delta }(1,0,0)f(z)\equiv D^{\delta
}f(z), $ $\left( \delta \in \mathrm{%
\mathbb{N}
}_{0}:=\mathrm{%
\mathbb{N}
}\cup \{0\}\right) $the Salagean differential operator [42].

\item $\mathrm{{\mathcal{J}}}_{1}^{\delta }(\lambda ,0,0)f(z)\equiv
D_{\lambda }^{\delta }f(z),$ $\left( \delta \in \mathrm{%
\mathbb{N}
}_{0}\right) $ the generalized Salagean differential operator introduced by
Al-Oboudi [2].

\item $\mathrm{{\mathcal{J}}}_{1}^{\delta }(\lambda ,\mu ,0)f(z)\equiv
D_{\lambda ,\mu }^{\delta }f(z),$ the operator studied by Deniz and Orhan
[18], in special case $0\leq \mu \leq \lambda \leq 1$ the operator was
studied firstly Raducanu and Orhan [39].

\item $\mathrm{{\mathcal{J}}}_{1}^{\delta }(1,0,l)f(z)\equiv I_{l}^{\delta
}f(z),$ $\left( \delta \in \mathrm{%
\mathbb{N}
}_{0}\right) $ the operator considered by Cho and Srivastava [15] and Cho
and Kim [16].

\item $\mathrm{{\mathcal{J}}}_{1}^{\delta }(1,0,1)f(z)\equiv I^{\delta
}f(z), $ $\left( \delta \in \mathrm{%
\mathbb{N}
}_{0}\right) $ the operator investigated by Uralegaddi and Somonatha [54].

\item $\mathrm{{\mathcal{J}}}_{1}^{\delta }(\lambda ,0,0)f(z)\equiv
D_{\lambda }^{\delta }f(z),$ $\left( \delta \in \mathrm{%
\mathbb{R}
}^{+}\cup \{0\}\right) $ the operator studied by Acu and Owa [1].

\item $\mathrm{{\mathcal{J}}}_{1}^{\delta }(\lambda ,0,l)f(z)\equiv I(\delta
,\lambda ,l)f(z),$ $\left( \delta \in \mathrm{%
\mathbb{R}
}^{+}\cup \{0\}\right) $ the operator introduced by Cata\c{s} [11].

\item $J_{p}^{\delta }(1,0,0)f(z)\equiv D_{p}^{\delta }f(z)$\textbf{$,$}$%
\left( \delta \in \mathrm{%
\mathbb{N}
}_{0}\right) $ the operator considered by Shenan \textit{et al. }[45].

\item $J_{p}^{\delta }(\lambda ,0,0)f(z)\equiv D_{\lambda ,p}^{\delta }f(z)$%
\textbf{$,$}$\left( \delta \in \mathrm{%
\mathbb{N}
}_{0}\right) $ the operator investigated by Kwon [25].

\item $\mathrm{{\mathcal{J}}}_{p}^{\delta }(1,0,l)f(z)\equiv I_{p}(\delta
,l)f(z),$ the operator considered by Kumar\textit{\ et al.} [48].

\item $\mathrm{{\mathcal{J}}}_{p}^{\delta }(\lambda ,0,l)f(z)\equiv
I_{p}(\delta ,\lambda ,l)f(z),$ the operator studied recently by Cata\c{s} 
\textit{et al}. [12].
\end{enumerate}

For special values of parameters $\lambda ,\mu ,l$ and $p$, from the
operator $\mathrm{{\mathcal{J}}}_{p}^{\delta }(\lambda ,\mu ,l)$ the
following new operators can be obtained:

$\bullet \;\mathrm{{\mathcal{J}}}_{p}^{\delta }(\lambda ,\mu ,1)\equiv 
\mathrm{{\mathcal{J}}}_{p}^{\delta }(\lambda ,\mu )$

$\bullet \;\mathrm{{\mathcal{J}}}_{1}^{\delta }(\lambda ,\mu ,l)\equiv 
\mathrm{{\mathcal{J}}}^{\delta }(\lambda ,\mu ,l).$

Now, by making use of the operator $\mathrm{{\mathcal{J}}}_{p}^{\delta
}(\lambda ,\mu ,l),$ we define a new subclass of functions belonging to the
class $\mathrm{{\mathcal{A}}}(n,p)$.

\noindent \textbf{Definition 1.2. }Let\textbf{\ }$\lambda \;\geq \mu \geq
0;\;l,\delta \geq 0;\;p\in \mathrm{%
\mathbb{N}
}$\textbf{\ }and for the parameters\textbf{\ $\sigma ,\;$}$A$ and $B$ such
that%
\begin{equation}
-1\leq A<B\leq 1,\;\;0<B\leq 1\text{ and }0\leq \sigma <p,  \tag{1.7}
\end{equation}

\noindent we say that a function $f(z)\in \mathrm{{\mathcal{A}}}(n,p)$ is in
the class $\mathcal{P}_{\lambda ,\mu ,l}^{\delta }(A,B;\sigma ,p)$ if it
satisfies the following subordination condition:%
\begin{equation}
\frac{1}{p-\sigma }\left( \frac{[\mathrm{{\mathcal{J}}}_{p}^{\delta
}(\lambda ,\mu ,l)f(z)]^{\prime }}{z^{p-1}}-\sigma \right) \prec \frac{1+Az}{%
1+Bz}\text{\ \ }(z\in \mathrm{{\mathcal{U}}}).  \tag{1.8}
\end{equation}%
If the following inequality holds true,%
\begin{equation}
\left\vert \frac{\frac{\lbrack \mathrm{{\mathcal{J}}}_{p}^{\delta }(\lambda
,\mu ,l)f(z)]^{\prime }}{z^{p-1}}-p}{B\frac{[\mathrm{{\mathcal{J}}}%
_{p}^{\delta }(\lambda ,\mu ,l)f(z)]^{\prime }}{z^{p-1}}-[pB+(A-B)(p-\sigma
)]}\right\vert <1\text{\ \ }(z\in \mathrm{{\mathcal{U}}})  \tag{1.9}
\end{equation}%
the inequality (1.9) is equivalent the subordination condition (1.8).

We note that by specializing the parameters $\lambda ,\mu ,l,\delta ,\sigma
,A,B$ and $p,$ the subclass $\mathcal{P}_{\lambda ,\mu ,l}^{\delta
}(A,B;\sigma ,p)$ reduces to several well-known subclasses of analytic
functions. These subclasses are:

\begin{enumerate}
\item $\mathcal{P}_{\lambda ,\mu ,l}^{0}(-1,1;0,1)\equiv \mathcal{P}%
_{0,0,l}^{1}(-1,1;0,1)=\mathrm{{\mathcal{R}}}$ (\textit{see Mac-Gregor}
[31]);

\item $\mathcal{P}_{\lambda ,\mu ,l}^{0}(A,B;\sigma ,p)\equiv \mathcal{P}%
_{0,0,l}^{1}(A,B;\sigma ,p)\equiv \mathrm{{\mathcal{S}}}_{p}(A,B,\sigma )$ (%
\textit{see Aouf }[3]);

\item $\mathcal{P}_{\lambda ,\mu ,l}^{0}(-1,1;\sigma ,p)\equiv \mathcal{P}%
_{0,0,l}^{1}(-1,1;\sigma ,p)\equiv \mathrm{{\mathcal{S}}}_{p}(\sigma )$ (%
\textit{see Owa }[35]);

\item $\mathcal{P}_{\lambda ,\mu ,l}^{0}(-1,1-\frac{1}{\alpha };0,1)\equiv 
\mathcal{P}_{0,0,l}^{1}(-1,1-\frac{1}{\alpha };0,1)\equiv \mathrm{{\mathcal{S%
}}}(\alpha )$ $\left( \alpha >\frac{1}{2}\right) $ (\textit{see Goel} [20]);

\item $\mathcal{P}_{\lambda ,\mu ,l}^{0}(-1,1-\frac{1}{\alpha };0,p)\equiv 
\mathcal{P}_{0,0,l}^{1}(-1,1-\frac{1}{\alpha };0,p)\equiv \mathrm{{\mathcal{S%
}}}_{p}(\alpha )$ $\left( \alpha >\frac{1}{2}\right) $ (\textit{see Sohi}
[49]);

\item $\mathcal{P}_{\lambda ,\mu ,l}^{0}(A,B;0,p)\equiv \mathcal{P}%
_{0,0,l}^{1}(A,B;0,p)\equiv \mathrm{{\mathcal{S}}}_{p}(A,B)$ (\textit{see
Chen} [13]);

\item $\mathcal{P}_{\lambda ,\mu ,l}^{0}(A,B;0,1)\equiv \mathcal{P}%
_{0,0,l}^{1}(A,B;0,1)\equiv \mathrm{{\mathcal{R}}}(A,B),\;\left( -1\leq
B<A\leq 1\right) $ (\textit{see Mehrok} [32]);

\item $\mathcal{P}_{\lambda ,\mu ,l}^{0}(-\gamma ,\gamma ;0,1)\equiv 
\mathcal{P}_{0,0,l}^{1}(-\gamma ,\gamma ;0,1)\equiv \mathrm{{\mathcal{R}}}%
_{(\gamma )},$ $\left( 0<\gamma \leq 1\right) $ (\textit{see Padmanabhan}
[38] \textit{and} \textit{Caplinger and Causey }[10]);

\item $\mathcal{P}_{\lambda ,\mu ,l}^{0}((2\beta -1)\gamma ,\gamma
;0,1)\equiv \mathcal{P}_{0,0,l}^{1}((2\beta -1)\gamma ,\gamma ;0,1)\equiv 
\mathrm{{\mathcal{R}}}_{\beta ,\gamma },$ $\left( 0\leq \beta <1,\;0<\gamma
\leq 1\right) $ (\textit{see Juneja and Mogra} [23]);

\item $\mathcal{P}_{\lambda ,\mu ,l}^{0}((2a-1)b,b;0,p)\equiv \mathcal{P}%
_{0,0,l}^{1}((2a-1)b,b;0,p)\equiv \mathrm{{\mathcal{S}}}_{p}(a,b),$ $\left(
0\leq a<1,\;0<b\leq 1\right) $ (\textit{see Owa }[36]);

\item $\mathcal{P}_{\lambda ,\mu ,l}^{0}((\gamma -1)\beta ,\alpha \beta
;0,1)\equiv \mathcal{P}_{0,0,l}^{1}((\gamma -1)\beta ,\alpha \beta
;0,1)\equiv \mathrm{{\mathcal{L}}}(\alpha ,\beta ,\gamma ),$ $\left( 0\leq
\alpha \leq 1,\;0<\beta \leq 1,0\leq \gamma <1\right) $ (\textit{see Kim and
Lee }[24]).
\end{enumerate}

Furthermore, we say that a function $f(z)\in \mathcal{P}_{\lambda ,\mu
,l}^{\delta }(A,B;\sigma ,p)$ is in the subclass $\widetilde{\mathcal{P}}%
_{\lambda ,\mu ,l}^{\delta }(A,B;\sigma ,p)$ if $f(z)$ is of the following
form:%
\begin{equation}
f(z)=z^{p}-\sum_{k=n+p}^{\infty }\left\vert a_{k}\right\vert z^{k}\text{\ \ }%
\left( p,n\in \mathrm{%
\mathbb{N}
}:=\{1,2,3,...\}\right) .  \tag{1.10}
\end{equation}%
Thus, by specializing the parameters $\lambda ,\mu ,l,\delta ,\sigma ,A,B$
and $p,$ we obtain the following familiar subclasses of analytic functions
in $\mathrm{{\mathcal{U}}}$ with \textit{negative }coefficients:

\begin{enumerate}
\item $\widetilde{\mathcal{P}}_{\lambda ,\mu ,l}^{0}(-1,1,;\alpha ,1)\equiv 
\widetilde{\mathcal{P}}_{0,0,l}^{1}(-1,1,;\alpha ,1)\equiv \mathcal{P}^{\ast
}(\alpha )$ $\left( 0\leq \alpha <1\right) $ (\textit{see for }$f\in \mathrm{%
{\mathcal{A}}}$ \textit{Sarangi and Uralegaddi }[43] \textit{and for $f\in 
\mathrm{{\mathcal{A}}}(n)$ Sekine and Owa }[44]);

\item $\widetilde{\mathcal{P}}_{\lambda ,\mu ,l}^{0}((\gamma -1)\beta
,\alpha \beta ;0,1)\equiv \widetilde{\mathcal{P}}_{0,0,l}^{1}((\gamma
-1)\beta ,\alpha \beta ;0,1)\equiv \mathrm{{\mathcal{L}}}^{\ast }(\alpha
,\beta ,\gamma )$ $\left( 0\leq \alpha \leq 1,\;0<\beta \leq 1,0\leq \gamma
<1\right) $ (\textit{see Kim and Lee} [24]);

\item $\widetilde{\mathcal{P}}_{\lambda ,\mu ,l}^{0}(A,B;\sigma ,p)\equiv 
\widetilde{\mathcal{P}}_{0,0,l}^{1}(A,B;\sigma ,p)\equiv \mathcal{P}^{\ast
}(p,A,B,\sigma )$ (\textit{see Aouf} [4]);

\item $\widetilde{\mathcal{P}}_{\lambda ,\mu ,l}^{0}(-\beta ,\beta
;0,1)\equiv \widetilde{\mathcal{P}}_{0,0,l}^{1}(-\beta ,\beta ;0,1)\equiv 
\mathrm{{\mathcal{D}}}^{\ast }(\beta )$ $0<\beta \leq 1$ (\textit{see Kim
and Lee} [24]);

\item $\widetilde{\mathcal{P}}_{\lambda ,\mu ,l}^{0}(-\beta ,\beta ;\sigma
,1)\equiv \widetilde{\mathcal{P}}_{0,0,l}^{1}(-\beta ,\beta ;\sigma
,1)\equiv \mathcal{P}^{\ast }(\sigma ,\beta )$ $\left( 0\leq \sigma
<p;\;0<\beta \leq 1\right) $ (\textit{see Gupta and Jain} [22]);

\item $\widetilde{\mathcal{P}}_{\lambda ,\mu ,l}^{0}(-\beta ,\beta ;\alpha
,p)\equiv \widetilde{\mathcal{P}}_{0,0,l}^{1}(-\beta ,\beta ;\alpha
,p)\equiv \mathcal{P}_{p}^{\ast }(\alpha ,\beta )$ $\left( 0\leq \alpha
<p;\;0<\beta \leq 1\right) $ (\textit{see Aouf } [6]);

\item $\widetilde{\mathcal{P}}_{\lambda ,\mu ,l}^{0}(A,B;0,p)\equiv 
\widetilde{\mathcal{P}}_{0,0,l}^{1}(A,B;0,p)\equiv \mathcal{P}^{\ast
}(p,A,B) $ (\textit{see} \textit{Shukla and Dashrath} [46]);

\item $\widetilde{\mathcal{P}}_{\lambda ,\mu ,l}^{0}(-1,1;\sigma ,p)\equiv 
\widetilde{\mathcal{P}}_{0,0,l}^{1}(-1,1;\sigma ,p)\equiv \mathrm{{\mathcal{F%
}}}_{p}(1,\beta )$ $\left( 0\leq \sigma <p;\;p\in \mathrm{%
\mathbb{N}
}\right) $ (\textit{see for} $f\in \mathrm{{\mathcal{A}}}$ \textit{Lee et al.%
} [26] \textit{and for $f\in \mathrm{{\mathcal{A}}}(n)$} \textit{Yaguchi et
al}.[55]).

\item $\widetilde{\mathcal{P}}_{\lambda ,\mu ,l}^{0}((2a-1)b,b;0,p)\equiv 
\widetilde{\mathcal{P}}_{0,0,l}^{1}((2a-1)b,b;0,p)\equiv \mathrm{{\mathcal{T}%
}}_{p}(a,b)$ $\left( 0\leq a<1,\;0<b\leq 1\right) $ (\textit{see Owa }[36]);
\end{enumerate}

In our present paper, we shall make use of the familiar \textit{integral
operator }$\mathrm{{\mathcal{I}}}_{\vartheta ,p}$ defined by (see, for
details, [9, 27, 30]; see also [54])%
\begin{equation}
(\mathrm{{\mathcal{I}}}_{\vartheta ,p})(z):=\frac{\vartheta +p}{z^{p}}%
\int_{0}^{z}t^{\vartheta -1}f(t)dt\text{\ \ \ }f\in \mathrm{{\mathcal{A}}}%
(n,p);\;\vartheta +p>0;\;p\in \mathrm{%
\mathbb{N}
})  \tag{1.11}
\end{equation}%
as well as the fractional calculus operator $\mathrm{{\mathcal{D}}}_{z}^{\nu
}$ for which it is well known that (see, for details, [37, 50] and [53]; see
also Section 7)%
\begin{equation}
\mathrm{{\mathcal{D}}}_{z}^{\nu }\{z^{p}\}=\frac{\Gamma (\rho +1)}{\Gamma
(\rho +1-\nu )}z^{\rho -\nu }\;\;(\rho >-1;\;\nu \in \mathrm{%
\mathbb{R}
})  \tag{1.12}
\end{equation}%
in terms of Gamma function.

The main object of the present paper is to investigate the various important
properties and characteristics of two subclasses of $\mathrm{{\mathcal{A}}}%
(n,p)$ of normalized analytic functions in $\mathrm{{\mathcal{U}}}$ with
negative and positive coefficients, which are introduced here by making use
of the multiplier transformations $\mathrm{{\mathcal{J}}}_{p}^{\delta
}(\lambda ,\mu ,l)$ defined by (1.4). Inclusion relationships for the class $%
\mathcal{P}_{\lambda ,\mu ,l}^{\delta }(A,B;\sigma ,p)$ are investigated by
applying the techniques of convolution. Furthermore, several properties
involving generalized neighborhoods and partial sums for functions belonging
to these subclasses are investigated. We also derive many results for the
Quasi- convolution of functions belonging to the class $\widetilde{\mathcal{P%
}}_{\lambda ,\mu ,l}^{\delta }(A,B;\sigma ,p)$\textbf{$.$} Finally, some
applications of fractional calculus operators are considered. Relevant
connections of the definitions and results presented here with those
obtained in several earlier works are also pointed out.

\noindent \textbf{Remark 1.2. }Throughout our present investigation, we
tacitly assume that the parametric constraints listed in (1.6), (1.7) and
Definition 1.1 are satisfied.

\section{INCLUSION PROPERTIES OF THE FUNCTION CLASS $\mathcal{P}_{\protect%
\lambda ,\protect\mu ,l}^{\protect\delta }(A,B;\protect\sigma ,p)$}

For proving our first inclusion result, we shall need the following lemmas.

\noindent \textbf{Lemma 2.1} (\textit{See Fejer} [19] \textit{or Ruscheweyh}
[41])\textbf{.} Assume $a_{1}=1$ and $a_{m}\geq 0$ for $m\geq 2,$ such that $%
\{a_{m}\}$ is a convex decreasing sequence, i.e., $a_{m}-2a_{m+1}+a_{m+2}%
\geq 0$ and $a_{m+1}-a_{m+2}\geq 0$ for $m\in \mathrm{%
\mathbb{N}
}.$ Then%
\begin{equation}
\Re \left\{ \sum_{m=1}^{\infty }a_{m}z^{m-1}\right\} >\frac{1}{2}  \tag{2.1}
\end{equation}%
for all $z\in \mathrm{{\mathcal{U}}}.$

\noindent \textbf{Lemma 2.2 }(\textit{See Liu} [28])\textbf{. }Let $-1\leq
A_{2}\leq A_{1}<B_{1}\leq B_{2}\leq 1.$ Then, we can write the following
subordination result:%
\begin{equation*}
\frac{1+B_{1}z}{1+A_{1}z}\prec \frac{1+B_{2}z}{1+A_{2}z}.
\end{equation*}

\noindent \textbf{Lemma 2.3. }If $\left[ \lambda -\mu \geq \frac{l+p}{2p}%
\mathrm{\;or\;\;}\lambda =\mu =0\right] ,$ then $\Re \left\{
1+\sum_{k=n+p}^{\infty }\frac{1}{\Phi _{p}^{k}(1,\lambda ,\mu ,l)}%
z^{k-p}\right\} >\frac{1}{2}$ for all $z\in \mathrm{{\mathcal{U}}}.$

\noindent \textbf{Proof. }Define:%
\begin{equation*}
q(z)=1+\sum_{k=2}^{\infty }\frac{1}{\Phi _{p}^{k+n+p-2}(1,\lambda ,\mu ,l)}%
z^{k+n-2}=1+\sum_{k=2}^{\infty }B_{k}z^{k+n-2}
\end{equation*}%
where%
\begin{equation}
B_{k}=\frac{1}{\Phi _{p}^{k+n+p-2}(1,\lambda ,\mu ,l)}=\frac{p+l}{%
(k+n-2)(\lambda \mu (k+n+p-2)+\lambda -\mu )+p+l}  \tag{2.2}
\end{equation}%
for all $n,p\in \mathrm{%
\mathbb{N}
},\;k\geq 2.$

\noindent Since the values $k,\;l,\;p,\;\lambda $ and $\mu $ are positive,
we have $B_{k}>0$ for all $k\in \mathrm{%
\mathbb{N}
}.$ We can easily find that%
\begin{equation}
B_{k+1}=\frac{p+l}{(k+n-1)(\lambda \mu (k+n+p-1)+\lambda -\mu )+p+l} 
\tag{2.3}
\end{equation}%
\begin{equation}
B_{k+2}=\frac{p+l}{(k+n)(\lambda \mu (k+n+p)+\lambda -\mu )+p+l}  \tag{2.4}
\end{equation}%
and thus from (2.3) and (2.4), we can see that

\begin{equation*}
B_{k+1}-B_{k+2}\geq 0
\end{equation*}%
for all $k\in \mathrm{%
\mathbb{N}
}.$ Next, we show that the inequality

\begin{equation}
B_{k}-2B_{k+1}+B_{k+2}\geq 0  \tag{2.5}
\end{equation}%
holds for all $k\in \mathrm{%
\mathbb{N}
}.$ Using (2.2), (2.3) and (2.4) we find that%
\begin{eqnarray*}
{B_{k}-2B_{k+1}+B_{k+2}} &=&{\frac{p+l}{(k+n-2)(\lambda \mu
(k+n+p-2)+\lambda -\mu )+p+l}} \\
&&-{2\frac{p+l}{(k+n-1)(\lambda \mu (k+n+p-1)+\lambda -\mu )+p+l}} \\
&&{+\frac{p+l}{(k+n)(\lambda \mu (k+n+p)+\lambda -\mu )+p+l}} \\
&=&\frac{2(\lambda \mu )^{2}\left[ 3(k+n)(k+n+p-2)+p^{2}-3p+2\right] }{%
C_{2}C_{1}C_{0}} \\
&&+\frac{2\lambda \mu \left[ 3(\lambda -\mu )(k+n-1)+p(2(\lambda -\mu )-1)-l%
\right] +(\lambda -\mu )^{2}}{C_{2}C_{1}C_{0}}
\end{eqnarray*}%
where $C_{i}=\left[ (k+n-i)(\lambda \mu (k+n+p-i)+\lambda -\mu )+p+l\right] $
and from the hypothesis of Lemma 2.3, we deduce that (2.5) holds for all $%
k\in \mathrm{%
\mathbb{N}
}.$ Thus the sequence $\{B_{k}\}$ is convex decreasing and by Lemma 2.1 we
obtain that%
\begin{equation*}
\Re \{q(z)\}=\Re \left\{ 1+\sum_{k=2}^{\infty }B_{k}z^{k+n-2}\right\} =\Re
\left\{ 1+\sum_{k=n+p}^{\infty }\frac{1}{\Phi _{p}^{k}(1,\lambda ,\mu ,l)}%
z^{k-p}\right\} >\frac{1}{2}
\end{equation*}%
for all $z\in \mathrm{{\mathcal{U}}}.$ The proof of Lemma 2.3 is completed.

\noindent \textbf{Theorem 2.1. }If $\left[ \lambda -\mu \geq \frac{l+p}{2p}%
\mathrm{\;or\;\;}\lambda =\mu =0\right] $ and $\delta \geq 0,$ then%
\begin{equation}
\mathcal{P}_{\lambda ,\mu ,l}^{\delta +1}(A,B;\sigma ,p)\subseteq \mathcal{P}%
_{\lambda ,\mu ,l}^{\delta }(A,B;\sigma ,p).  \tag{2.6}
\end{equation}

\noindent \textbf{Proof.} Let $f\in \mathcal{P}_{\lambda ,\mu ,l}^{\delta
+1}(A,B;\sigma ,p).$ Using the definition of $\mathcal{P}_{\lambda ,\mu
,l}^{\delta +1}(A,B;\sigma ,p)$ we obtain that%
\begin{equation}
\frac{1}{p-\sigma }\left( \frac{z[\mathrm{{\mathcal{J}}}_{p}^{\delta
+1}(\lambda ,\mu ,l)f(z)]^{\prime }}{z^{p}}-\sigma \right) \prec \frac{1+Az}{%
1+Bz}\text{ \ }(z\in \mathrm{{\mathcal{U}}}).  \tag{2.7}
\end{equation}%
Applying definition of $\mathrm{{\mathcal{J}}}_{p}^{\delta }(\lambda ,\mu
,l)f(z)$ and the properties of convolution we find that%
\begin{eqnarray*}
&&\frac{z[\mathrm{{\mathcal{J}}}_{p}^{\delta }(\lambda ,\mu ,l)f(z)]^{\prime
}}{pz^{p}} \\
&=&\left( 1+\sum_{k=n+p}^{\infty }\frac{1}{\Phi _{p}^{k}(1,\lambda ,\mu ,l)}%
z^{k-p}\right) \ast \left( 1+\sum_{k=n+p}^{\infty }\frac{k}{p}\Phi
_{p}^{k}(\delta +1,\lambda ,\mu ,l)a_{k}z^{k-p}\right) \\
&=&\left( 1+\sum_{k=n+p}^{\infty }\frac{1}{\Phi _{p}^{k}(1,\lambda ,\mu ,l)}%
z^{k-p}\right) \ast \left( \frac{z[\mathrm{{\mathcal{J}}}_{p}^{\delta
+1}(\lambda ,\mu ,l)f(z)]^{\prime }}{pz^{p}}\right) .
\end{eqnarray*}%
Therefore from the last equalities and (1.8) we get%
\begin{eqnarray}
&&\frac{1}{p-\sigma }\left( \frac{[\mathrm{{\mathcal{J}}}_{p}^{\delta
}(\lambda ,\mu ,l)f(z)]^{\prime }}{z^{p-1}}-\sigma \right)  \TCItag{2.8} \\
&=&\frac{1}{p-\sigma }\left\{ p\left[ \left( 1+\sum_{k=n+p}^{\infty }\frac{1%
}{\Phi _{p}^{k}(1,\lambda ,\mu ,l)}z^{k-p}\right) \ast \left( \frac{z[%
\mathrm{{\mathcal{J}}}_{p}^{\delta +1}(\lambda ,\mu ,l)f(z)]^{\prime }}{%
pz^{p}}\right) \right] -\sigma \right\}  \notag \\
&=&\left( 1+\sum_{k=n+p}^{\infty }\frac{1}{\Phi _{p}^{k}(1,\lambda ,\mu ,l)}%
z^{k-p}\right) \ast \frac{1}{p-\sigma }\left( \frac{z[\mathrm{{\mathcal{J}}}%
_{p}^{\delta +1}(\lambda ,\mu ,l)f(z)]^{\prime }}{z^{p}}-\sigma \right) 
\notag \\
&=&q(z)\ast \frac{1+Aw(z)}{1+Bw(z)}.  \notag
\end{eqnarray}%
where $\left\vert w(z)\right\vert <1$ and $w(0)=0.$ From the Herglotz
theorem and Lemma 2.3 we thus obtain%
\begin{equation*}
q(z)=\int_{\left\vert x\right\vert =1}\frac{d\varpi (x)}{1-xz}\;\;(z\in 
\mathrm{{\mathcal{U}}}),
\end{equation*}%
when $\varpi (x)$ is a probability measure on the unit circle $\left\vert
x\right\vert =1,$ that is%
\begin{equation*}
\int_{\left\vert x\right\vert =1}d\varpi (x)=1.
\end{equation*}%
It follows from (2.8) that%
\begin{equation*}
\frac{1}{p-\sigma }\left( \frac{z[\mathrm{{\mathcal{J}}}_{p}^{\delta
}(\lambda ,\mu ,l)f(z)]^{\prime }}{z^{p}}-\sigma \right) =\int_{\left\vert
x\right\vert =1}\frac{1+Axz}{1+Bxz}d\varpi (x)\prec \frac{1+Az}{1+Bz}
\end{equation*}%
because $\frac{1+Az}{1+Bz}$ is convex univalent in $\mathrm{{\mathcal{U}}}.$
Hence we conclude that%
\begin{equation*}
\mathcal{P}_{\lambda ,\mu ,l}^{\delta +1}(A,B;\sigma ,p)\subseteq \mathcal{P}%
_{\lambda ,\mu ,l}^{\delta }(A,B;\sigma ,p),
\end{equation*}%
which completes the proof of Theorem 2.1.

\noindent \textbf{Theorem 2.2. }If $\delta \geq 0$ and $-1\leq A_{2}\leq
A_{1}<B_{1}\leq B_{2}\leq 1,$ then%
\begin{equation}
\mathcal{P}_{\lambda ,\mu ,l}^{\delta +1}(B_{1},A_{1};\sigma ,p)\subseteq 
\mathcal{P}_{\lambda ,\mu ,l}^{\delta }(B_{2},A_{2};\sigma ,p).  \tag{2.9}
\end{equation}

\noindent \textbf{Proof. }Making use of Lemma 2.2, we can write%
\begin{equation*}
\mathcal{P}_{\lambda ,\mu ,l}^{\delta }(B_{1},A_{1};\sigma ,p)\subseteq 
\mathcal{P}_{\lambda ,\mu ,l}^{\delta }(B_{2},A_{2};\sigma ,p).
\end{equation*}%
Using (2.6) and (2.9), we have%
\begin{equation*}
\mathcal{P}_{\lambda ,\mu ,l}^{\delta +1}(B_{1},A_{1};\sigma ,p)\subseteq 
\mathcal{P}_{\lambda ,\mu ,l}^{\delta }(B_{1},A_{1};\sigma ,p)\subseteq 
\mathcal{P}_{\lambda ,\mu ,l}^{\delta }(B_{2},A_{2};\sigma ,p)
\end{equation*}%
so, we obtain%
\begin{equation*}
\mathcal{P}_{\lambda ,\mu ,l}^{\delta +1}(B_{1},A_{1};\sigma ,p)\subseteq 
\mathcal{P}_{\lambda ,\mu ,l}^{\delta }(B_{2},A_{2};\sigma ,p).
\end{equation*}%
Thus, the proof is complete.

\section{BASIC PROPERTIES OF THE FUNCTION CLASS\textbf{\ }$\protect%
\widetilde{\mathcal{P}}_{\protect\lambda ,\protect\mu ,l}^{\protect\delta %
}(A,B;\protect\sigma ,p)$}

We first determine a necessary and sufficient condition for a function $%
f(z)\in \mathrm{{\mathcal{A}}}(n,p)$ of the form (1.10) to be in the class $%
\widetilde{\mathcal{P}}_{\lambda ,\mu ,l}^{\delta }(A,B;\sigma ,p).$

\noindent \textbf{Theorem 3.1. }Let the function $f(z)\in \mathrm{{\mathcal{A%
}}}(n,p)$ be defined by (1.10). Then the function $f(z)$ is in the class $%
\widetilde{\mathcal{P}}_{\lambda ,\mu ,l}^{\delta }(A,B;\sigma ,p)$\textbf{\ 
}if and only if%
\begin{equation}
\sum_{k=n+p}^{\infty }k(1+B)\Phi _{p}^{k}(\delta ,\lambda ,\mu ,l)\left\vert
a_{k}\right\vert \leq (B-A)(p-\sigma )  \tag{3.1}
\end{equation}%
where $\Phi _{p}^{k}(\delta ,\lambda ,\mu ,l)$ is given by (1.6).

\noindent \textbf{Proof. }If the condition (3.1) hold true, we find from
(1.10) and (3.1) that%
\begin{eqnarray*}
&&\left\vert [\mathrm{{\mathcal{J}}}_{p}^{\delta }(\lambda ,\mu
,l)f(z)]^{\prime }-pz^{p-1}\right\vert -\left\vert B[\mathrm{{\mathcal{J}}}%
_{p}^{\delta }(\lambda ,\mu ,l)f(z)]^{\prime }-z^{p-1}[pB+(A-B)(p-\sigma
)]\right\vert \\
&=&\left\vert -\sum_{k=n+p}^{\infty }k\Phi _{p}^{k}(\delta ,\lambda ,\mu
,l)\left\vert a_{k}\right\vert z^{k-1}\right\vert -\left\vert (B-A)(p-\sigma
)z^{p-1}-B\sum_{k=n+p}^{\infty }k\Phi _{p}^{k}(\delta ,\lambda ,\mu
,l)\left\vert a_{k}\right\vert z^{k-1}\right\vert \\
&\leq &\sum_{k=n+p}^{\infty }k(1+B)\Phi _{p}^{k}(\delta ,\lambda ,\mu
,l)\left\vert a_{k}\right\vert -(B-A)(p-\sigma )\leq 0\text{\ \ }\left( z\in
\partial \mathrm{{\mathcal{U}}}=\{z:\;z\in \mathrm{%
\mathbb{C}
}\;\mathrm{and\;}\left\vert z\right\vert =1\}\right) .
\end{eqnarray*}%
Hence, by the \textit{Maximum Modulus Theorem}, we have%
\begin{equation*}
f(z)\in \widetilde{\mathcal{P}}_{\lambda ,\mu ,l}^{\delta }(A,B;\sigma ,p).
\end{equation*}%
Conversely, assume that the function $f(z)$ defined by (1.10) is in the
class $\widetilde{\mathcal{P}}_{\lambda ,\mu ,l}^{\delta }(A,B;\sigma ,p)$%
\textbf{$.$ }Then\textbf{\ }we have%
\begin{eqnarray}
&&\left\vert \frac{\frac{[\mathrm{{\mathcal{J}}}_{p}^{\delta }(\lambda ,\mu
,l)f(z)]^{\prime }}{z^{p-1}}-p}{B\frac{[\mathrm{{\mathcal{J}}}_{p}^{\delta
}(\lambda ,\mu ,l)f(z)]^{\prime }}{z^{p-1}}-[pB+(A-B)(p-\sigma )]}\right\vert
\TCItag{3.2} \\
&=&\left\vert \frac{\sum_{k=n+p}^{\infty }k\Phi _{p}^{k}(\delta ,\lambda
,\mu ,l)\left\vert a_{k}\right\vert z^{k-p}}{(B-A)(p-\sigma
)z^{p-1}-B\sum_{k=n+p}^{\infty }k\Phi _{p}^{k}(\delta ,\lambda ,\mu
,l)\left\vert a_{k}\right\vert z^{k-p}}\right\vert <1\;\;(z\in \mathrm{{%
\mathcal{U}}}).  \notag
\end{eqnarray}%
Now, since $\left\vert \Re (z)\right\vert \leq \left\vert z\right\vert $ for
all $z,$ we have%
\begin{equation}
\Re \left( \frac{\sum_{k=n+p}^{\infty }k\Phi _{p}^{k}(\delta ,\lambda ,\mu
,l)\left\vert a_{k}\right\vert z^{k-p}}{(B-A)(p-\sigma
)z^{p-1}-B\sum_{k=n+p}^{\infty }k\Phi _{p}^{k}(\delta ,\lambda ,\mu
,l)\left\vert a_{k}\right\vert z^{k-p}}\right) <1.  \tag{3.3}
\end{equation}%
We choose values of $z$ on the real axis so that the following expression:%
\begin{equation*}
\frac{\lbrack \mathrm{{\mathcal{J}}}_{p}^{\delta }(\lambda ,\mu
,l)f(z)]^{\prime }}{z^{p-1}}
\end{equation*}%
is real. Then, upon clearing the denominator in (3.3) and letting $%
z\rightarrow 1^{-}$ though real values, we get the following inequality%
\begin{equation*}
\sum_{k=n+p}^{\infty }k(1+B)\Phi _{p}^{k}(\delta ,\lambda ,\mu ,l)\left\vert
a_{k}\right\vert \leq (B-A)(p-\sigma ).
\end{equation*}%
This completes the proof of Theorem 3.1.

\noindent \textbf{Remark 3.1. }Since $\widetilde{\mathcal{P}}_{\lambda ,\mu
,l}^{\delta }(A,B;\sigma ,p)$ is contained in the function class $\mathcal{P}%
_{\lambda ,\mu ,l}^{\delta }(A,B;\sigma ,p)$\textbf{$,$ }a sufficient
condition for $f(z)$ defined by (1.1) to be in the class $\mathcal{P}%
_{\lambda ,\mu ,l}^{\delta }(A,B;\sigma ,p)$ is that it satisfies the
condition (3.1) of Theorem 3.1.

\noindent \textbf{Corollary 3.1. }Let the function $f(z)\in \mathrm{{%
\mathcal{A}}}(n,p)$ be defined by (1.10). If the function $f(z)\in 
\widetilde{\mathcal{P}}_{\lambda ,\mu ,l}^{\delta }(A,B;\sigma ,p),$ then%
\begin{equation}
\left\vert a_{k}\right\vert \leq \frac{(B-A)(p-\sigma )}{k(1+B)\Phi
_{p}^{k}(\delta ,\lambda ,\mu ,l)}\text{\ \ }\left( k,p\in \mathrm{%
\mathbb{N}
}\right) .  \tag{3.4}
\end{equation}%
The result is sharp for the function $f(z)$ given by%
\begin{equation}
f(z)=z^{p}-\frac{(B-A)(p-\sigma )}{k(1+B)\Phi _{p}^{k}(\delta ,\lambda ,\mu
,l)}z^{k}\text{\ \ }\left( k,p\in \mathrm{%
\mathbb{N}
}\right) .  \tag{3.5}
\end{equation}%
We next prove the following growth and distortion properties for the class $%
\widetilde{\mathcal{P}}_{\lambda ,\mu ,l}^{\delta }(A,B;\sigma ,p)$\textbf{$%
. $ }

\noindent \textbf{Remark 3.2.}

\begin{enumerate}
\item Putting $A=-\beta ,\;B=\beta ,\;\delta =0$ and $\sigma =\alpha $ in
Theorem 3.1, we obtain the corresponding result given earlier by Aouf [6].

\item Putting $A=(\gamma -1)\beta ,\;B=\alpha \beta ,\;\delta =0,\;p=n=1$
and $\sigma =0$ in Theorem 3.1, we obtain result of Kim and Lee [24].

\item Putting $\delta =0$ in Theorem 3.1, we obtain Theorem 1 in [4].

\item Putting $A=(2a-1)b,\;B=b,\;\delta =0,\;n=1$ and $\sigma =0$ in Theorem
3.1, we arrive at the Theorem of Owa [36].

\item Putting $\delta =0$ and $\sigma =0$ in Theorem 3.1, we obtain the
corresponding result due to Shukla and Dashrath [46].
\end{enumerate}

\noindent \textbf{Theorem 3.2. }If a function $f(z)$ be defined by (1.10) is
in the class $\widetilde{\mathcal{P}}_{\lambda ,\mu ,l}^{\delta }(A,B;\sigma
,p),$ then%
\begin{eqnarray}
&&~~~~~~\left( \frac{p!}{(p-q)!}-\frac{(B-A)(p-\sigma )(n+p-1)!}{(1+B)\Phi
_{p}^{n+p}(\delta ,\lambda ,\mu ,l)(n+p-q)!}\left\vert z\right\vert
^{n}\right) \left\vert z\right\vert ^{p-q}  \TCItag{3.6} \\
&\leq &\left\vert f^{(q)}(z)\right\vert \leq \left( \frac{p!}{(p-q)!}+\frac{%
(B-A)(p-\sigma )(n+p-1)!}{(1+B)\Phi _{p}^{n+p}(\delta ,\lambda ,\mu
,l)(n+p-q)!}\left\vert z\right\vert ^{n}\right) \left\vert z\right\vert
^{p-q}  \notag
\end{eqnarray}%
for $q\in \mathrm{%
\mathbb{N}
}_{0},\;p>q$ and all $z\in \mathrm{{\mathcal{U}}}.$ The result is sharp for
the function $f(z)$ given by%
\begin{equation}
f(z)=z^{p}-\frac{(B-A)(p-\sigma )}{(n+p)(1+B)\Phi _{p}^{n+p}(\delta ,\lambda
,\mu ,l)}z^{n+p}\text{\ \ }\left( p\in \mathrm{%
\mathbb{N}
}\right) .  \tag{3.7}
\end{equation}

\noindent \textbf{Proof. }In view of Theorem 3.1, we have%
\begin{equation*}
\frac{(n+p)(1+B)\Phi _{p}^{n+p}(\delta ,\lambda ,\mu ,l)}{(B-A)(p-\sigma
)(n+p)!}\sum_{k=n+p}^{\infty }k!\left\vert a_{k}\right\vert \leq
\sum_{k=n+p}^{\infty }\frac{k(1+B)\Phi _{p}^{k}(\delta ,\lambda ,\mu ,l)}{%
(B-A)(p-\sigma )}\left\vert a_{k}\right\vert \leq 1,
\end{equation*}%
which readily yields%
\begin{equation}
\sum_{k=n+p}^{\infty }k!\left\vert a_{k}\right\vert \leq \frac{%
(B-A)(p-\sigma )(n+p-1)!}{(1+B)\Phi _{p}^{n+p}(\delta ,\lambda ,\mu ,l)}%
\text{\ \ }\left( k,p\in \mathrm{%
\mathbb{N}
}\right) .  \tag{3.8}
\end{equation}%
Now, by differentiating both sides of (1.10) $q-$times with respect to $z,$
we obtain%
\begin{equation}
f^{(q)}(z)=\frac{p!}{(p-q)!}z^{p-q}-\sum_{k=n+p}^{\infty }\frac{k!}{(k-q)!}%
a_{k}z^{k-q}\;\;\left( q\in \mathrm{%
\mathbb{N}
}_{0};\;p>q\right) .  \tag{3.9}
\end{equation}%
Theorem 3.2 follows readily from (3.8) and (3.9).

\noindent Finally, it is easy to see that the bounds in (3.6) are attained
for the function $f(z)$ given by (3.7).

\section{\noindent INCLUSION RELATIONS INVOLVING NEIGHBORHOODS}

Following the earlier investigations (based upon the familiar concept of
neighborhoods of analytic functions) by Goodman [21], Ruscheweyh [40] and
others including Srivastava \textit{et al.} [50, 52], Orhan [33, 34], Deniz 
\textit{et al. }[17], Aouf \textit{et al}. [8] (see also [11]).

Firstly, we define the $(n,\eta )-$neighborhood of function $f(z)\in \mathrm{%
{\mathcal{A}}}(n,p)$ of the form (1.1) by means of Definition 4.1 below.

\noindent \textbf{Definition 4.1. }For $\eta >0$ and a non-negative sequence 
$\mathrm{{\mathcal{S}}}=\{s_{k}\}_{k=1}^{\infty },$ where%
\begin{equation}
s_{k}:=\frac{k(1+B)\Phi _{p}^{k}(\delta ,\lambda ,\mu ,l)}{(B-A)(p-\sigma )}%
\;\;(k\in \mathrm{%
\mathbb{N}
}).  \tag{4.1}
\end{equation}%
The $(n,\eta )-$neighborhood of a function $f(z)\in \mathrm{{\mathcal{A}}}%
(n,p)$ of the form (1.1) is defined as follows:%
\begin{equation}
\mathrm{{\mathcal{N}}}_{n,p}^{\eta }(f):=\left\{
g:\;g(z)=z^{p}+\sum_{k=n+p}^{\infty }b_{k}z^{k}\in \mathrm{{\mathcal{A}}}%
(n,p)\;\mathrm{and\;}\sum_{k=n+p}^{\infty }s_{k}\left\vert
b_{k}-a_{k}\right\vert \leq \eta \;(\eta >0)\right\} .  \tag{4.2}
\end{equation}

For $s_{k}=k,$ Definition 4.1 would correspond to the $\mathrm{{\mathcal{N}}}%
_{\eta }-$neighborhood considered by Ruscheweyh [40].

\noindent Our first result based upon the familiar concept of neighborhood
defined by (4.2).

\noindent \textbf{Theorem 4.1. }Let $f(z)\in \mathcal{P}_{\lambda ,\mu
,l}^{\delta }(A,B;\sigma ,p)$ be given by (1.1). If $f$ satisfies the
inclusion condition:%
\begin{equation}
\left( f(z)+\varepsilon z^{p}\right) \left( 1+\varepsilon \right) ^{-1}\in 
\mathcal{P}_{\lambda ,\mu ,l}^{\delta }(A,B;\sigma ,p)\text{\ \ }\left(
\varepsilon \in \mathrm{%
\mathbb{C}
};\;\left\vert \varepsilon \right\vert <\eta ;\;\eta >0\right) ,  \tag{4.3}
\end{equation}%
then%
\begin{equation}
\mathrm{{\mathcal{N}}}_{n,p}^{\eta }(f)\subset \mathcal{P}_{\lambda ,\mu
,l}^{\delta }(A,B;\sigma ,p).  \tag{4.4}
\end{equation}

\noindent \textbf{Proof. }It is not difficult to see that a function $f$ \
belongs to $\mathcal{P}_{\lambda ,\mu ,l}^{\delta }(A,B;\sigma ,p)$ if and
only if%
\begin{equation}
\frac{\lbrack \mathrm{{\mathcal{J}}}_{p}^{\delta }(\lambda ,\mu
,l)f(z)]^{\prime }-pz^{p-1}}{B[\mathrm{{\mathcal{J}}}_{p}^{\delta }(\lambda
,\mu ,l)f(z)]^{\prime }-z^{p-1}[pB+(A-B)(p-\sigma )]}\neq \tau \;\;\left(
z\in \mathrm{{\mathcal{U}}};\;\tau \in \mathrm{%
\mathbb{C}
},\;\left\vert \tau \right\vert =1\right) ,  \tag{4.5}
\end{equation}%
which is equivalent to%
\begin{equation}
{(f\ast h)(z)\diagup z^{p}}\neq 0\;\;(z\in \mathrm{{\mathcal{U}}}), 
\tag{4.6}
\end{equation}%
where for convenience,%
\begin{equation}
h(z):=z^{p}+\sum_{k=n+p}^{\infty }c_{k}z^{k}=z^{p}+\sum_{k=n+p}^{\infty }%
\frac{k(1+\tau B)\Phi _{p}^{k}(\delta ,\lambda ,\mu ,l)}{\tau (B-A)(p-\sigma
)}z^{k}.  \tag{4.7}
\end{equation}%
We easily find from (4.7) that%
\begin{equation}
\left\vert c_{k}\right\vert \leq \left\vert \frac{k(1+\tau B)\Phi
_{p}^{k}(\delta ,\lambda ,\mu ,l)}{\tau (B-A)(p-\sigma )}\right\vert \leq 
\frac{k(1+B)\Phi _{p}^{k}(\delta ,\lambda ,\mu ,l)}{(B-A)(p-\sigma )}%
\;\;(k\in \mathrm{%
\mathbb{N}
}).  \tag{4.8}
\end{equation}%
Furthermore, under the hypotheses of theorem, (4.3) and (4.6) yields the
following inequalities:%
\begin{equation*}
\frac{\left( (f(z)+\varepsilon z^{p})(1+\varepsilon )^{-1}\right) \ast h(z)}{%
z^{p}}\neq 0\;\;(z\in \mathrm{{\mathcal{U}}})
\end{equation*}%
or%
\begin{equation*}
\frac{f(z)\ast h(z)}{z^{p}}\neq \varepsilon \;\;(z\in \mathrm{{\mathcal{U}}}%
),
\end{equation*}%
which is equivalent to the following:%
\begin{equation}
\frac{f(z)\ast h(z)}{z^{p}}\geq \eta \;\;(z\in \mathrm{{\mathcal{U}}};\;\eta
>0).  \tag{4.9}
\end{equation}%
Now, if we let%
\begin{equation*}
g(z):=z^{p}+\sum_{k=n+p}^{\infty }b_{k}z^{k}\in \mathrm{{\mathcal{N}}}%
_{n,p}^{\eta }(f),
\end{equation*}%
then we have%
\begin{equation*}
{\left\vert \frac{\left( f(z)-g(z)\right) \ast h(z)}{z^{p}}\right\vert
=\left\vert \sum_{k=n+p}^{\infty }(a_{k}-b_{k})c_{k}z^{k-p}\right\vert }
\end{equation*}%
\begin{equation*}
{\;\leq \sum_{k=n+p}^{\infty }\frac{k(1+B)\Phi _{p}^{k}(\delta ,\lambda ,\mu
,l)}{(B-A)(p-\sigma )}\left\vert a_{k}-b_{k}\right\vert \left\vert
z\right\vert ^{k-p}<\eta }\text{\ \ }{(z\in \mathrm{{\mathcal{U}}};\;\eta
>0).}
\end{equation*}%
Thus, for any complex number $\tau $ such that $\left\vert \tau \right\vert
=1,$ we have 
\begin{equation*}
{(g\ast h)(z)\diagup z^{p}}\neq 0\;\;(z\in \mathrm{{\mathcal{U}}}),
\end{equation*}

\noindent which implies that $g\in \mathcal{P}_{\lambda ,\mu ,l}^{\delta
}(A,B;\sigma ,p).$ The proof is complete.

We now define the $(n,\eta )-$neighborhood of a function $f(z)\in \mathrm{{%
\mathcal{A}}}(n,p)$ of the form (1.10) as follows:

\noindent \textbf{Definition 4.2. }For $\eta >0,$ the $(n,\eta )-$%
neighborhood of a function $f(z)\in \mathrm{{\mathcal{A}}}(n,p)$ of the form
(1.10) is given by%
\begin{eqnarray}
{\widetilde{\mathrm{{\mathcal{N}}}}_{n,p}^{\eta }(f)}{:=} &&\left\{ {%
g:\;g(z)=z^{p}-\sum_{k=n+p}^{\infty }b_{k}z^{k}\in \mathrm{{\mathcal{A}}}%
(n,p)\;}\text{and}\right.  \TCItag{4.10} \\
&&\left. \text{ }\sum_{k=n+p}^{\infty }\frac{k(1+B)\Phi _{p}^{k}(\delta
,\lambda ,\mu ,l)}{(B-A)(p-\sigma )}\left\vert \left\vert b_{k}\right\vert
-\left\vert a_{k}\right\vert \right\vert \leq \eta \;\;(\eta >0)\right\} . 
\notag
\end{eqnarray}

Next, we prove

\noindent \textbf{Theorem 4.2. }If the function $f(z)$ defined by (1.10) is
in the class $\widetilde{\mathcal{P}}_{\lambda ,\mu ,l}^{\delta
+1}(A,B;\sigma ,p),$ then%
\begin{equation}
\widetilde{\mathrm{{\mathcal{N}}}}_{n,p}^{\eta }(f)\subset \widetilde{%
\mathcal{P}}_{\lambda ,\mu ,l}^{\delta }(A,B;\sigma ,p)  \tag{4.11}
\end{equation}%
where%
\begin{equation*}
\eta :=\frac{n[\lambda \mu (n+p)+\lambda -\mu ]}{n[\lambda \mu (n+p)+\lambda
-\mu ]+p+l}.
\end{equation*}%
The result is the best possible in the sense that $\eta $ cannot be
increased.

\noindent \textbf{Proof. }For a function $f(z)\in \widetilde{\mathcal{P}}%
_{\lambda ,\mu ,l}^{\delta +1}(A,B;\sigma ,p)$ of the form (1.10) Theorem
3.1 immediately yields%
\begin{equation}
\sum_{k=n+p}^{\infty }\frac{k(1+B)\Phi _{p}^{k}(\delta ,\lambda ,\mu ,l)}{%
(B-A)(p-\sigma )}\left\vert a_{k}\right\vert \leq \frac{p+l}{n[\lambda \mu
(n+p)+\lambda -\mu ]+p+l}.  \tag{4.12}
\end{equation}%
Similarly, by taking%
\begin{equation*}
g(z):=z^{p}-\sum_{k=n+p}^{\infty }\left\vert b_{k}\right\vert z^{k}\in 
\widetilde{\mathrm{{\mathcal{N}}}}_{n,p}^{\eta }(f)\;\;\left( \eta =\frac{%
n[\lambda \mu (n+p)+\lambda -\mu ]}{n[\lambda \mu (n+p)+\lambda -\mu ]+p+l}%
\right) ,
\end{equation*}%
we find from the definition (4.10) that%
\begin{equation}
\sum_{k=n+p}^{\infty }\frac{k(1+B)\Phi _{p}^{k}(\delta ,\lambda ,\mu ,l)}{%
(B-A)(p-\sigma )}\left\vert \left\vert b_{k}\right\vert -\left\vert
a_{k}\right\vert \right\vert \leq \eta \;\;(\eta >0).  \tag{4.13}
\end{equation}%
With the help of (4.12) and (4.13), we have%
\begin{eqnarray*}
{\sum_{k=n+p}^{\infty }\frac{k(1+B)\Phi _{p}^{k}(\delta ,\lambda ,\mu ,l)}{%
(B-A)(p-\sigma )}\left\vert b_{k}\right\vert } &\leq &{\sum_{k=n+p}^{\infty }%
\frac{k(1+B)\Phi _{p}^{k}(\delta ,\lambda ,\mu ,l)}{(B-A)(p-\sigma )}%
\left\vert b_{k}\right\vert } \\
&&{+\sum_{k=n+p}^{\infty }\frac{k(1+B)\Phi _{p}^{k}(\delta ,\lambda ,\mu ,l)%
}{(B-A)(p-\sigma )}\left\vert \left\vert b_{k}\right\vert -\left\vert
a_{k}\right\vert \right\vert } \\
&&{+\sum_{k=n+p}^{\infty }\frac{k(1+B)\Phi _{p}^{k}(\delta ,\lambda ,\mu ,l)%
}{(B-A)(p-\sigma )}\left\vert \left\vert b_{k}\right\vert -\left\vert
a_{k}\right\vert \right\vert } \\
&{\leq }&{\frac{p+l}{n[\lambda \mu (n+p)+\lambda -\mu ]+p+l}+\eta =1.}
\end{eqnarray*}%
Hence, in view of the Theorem 3.1 again, we see that $g(z)\in \widetilde{%
\mathcal{P}}_{\lambda ,\mu ,l}^{\delta +1}(A,B;\sigma ,p).$

To show the sharpness of the assertion of Theorem 4.2, we consider the
functions $f(z)$ and $g(z)$ given by%
\begin{equation}
f(z)=z^{p}-\left[ \frac{(B-A)(p-\sigma )}{(n+p)(1+B)\Phi _{p}^{n+p}(\delta
+1,\lambda ,\mu ,l)}\right] z^{n+p}\in \widetilde{\mathcal{P}}_{\lambda ,\mu
,l}^{\delta +1}(A,B;\sigma ,p)  \tag{4.14}
\end{equation}%
and%
\begin{equation}
g(z)=z^{p}-\left[ \frac{(B-A)(p-\sigma )}{(n+p)(1+B)\Phi _{p}^{n+p}(\delta
+1,\lambda ,\mu ,l)}+\frac{(B-A)(p-\sigma )}{(n+p)(1+B)\Phi
_{p}^{n+p}(\delta ,\lambda ,\mu ,l)}\eta ^{\ast }\right] z^{n+p}  \tag{4.15}
\end{equation}%
where $\eta ^{\ast }>\eta .$

\noindent Clearly, the function $g(z)$ belong to $\widetilde{\mathrm{{%
\mathcal{N}}}}_{n,p}^{\eta ^{\ast }}(f).$ On the other hand, we find from
Theorem 3.1 that $g(z)\notin \widetilde{\mathcal{P}}_{\lambda ,\mu
,l}^{\delta }(A,B;\sigma ,p).$This evidently completes the proof of Theorem
4.2.

\section{\noindent PARTIAL SUMS OF THE FUNCTION CLASS $\protect\widetilde{%
\mathcal{P}}_{\protect\lambda ,\protect\mu ,l}^{\protect\delta }(A,B;\protect%
\sigma ,p)$}

Following the earlier work by Silverman [47] and recently Liu [29] and Deniz 
\textit{et al}. [17], in this section we investigate the ratio of real parts
of functions involving (1.10) and its sequence of partial sums defined by%
\begin{equation}
\kappa _{m}(z)=\left\{ 
\begin{array}{l}
{z^{p},\;\;\;\;\;\;\;\;\;\;\;\;\;\;\;\;\;\;\;\;\ \ \ \ \ \ \
\;m=1,2,...,n+p-1;} \\ 
{z^{p}-\sum_{k=n+p}^{m}\left\vert a_{k}\right\vert
z^{k},\;\;\;m=n+p,n+p+1,....}%
\end{array}%
\right. \;\;(k\geq n+p;\;n,p\in \mathrm{%
\mathbb{N}
})  \tag{5.1}
\end{equation}%
and determine sharp lower bounds for $\Re \left\{ {f(z)\diagup \kappa _{m}(z)%
}\right\} ,\;\Re \left\{ {\kappa _{m}(z)\diagup f(z)}\right\} .$

\noindent \textbf{Theorem 5.1. }Let $f\in \mathrm{{\mathcal{A}}}(n,p)$ and $%
\kappa _{m}(z)$ be given by (1.10) and (5.1), respectively. Suppose also that%
\begin{equation}
\sum_{k=n+p}^{\infty }\theta _{k}\left\vert a_{k}\right\vert \leq 1~~\left( 
\mathrm{where\;}{\theta }_{k}=\frac{k(1+B)\Phi _{p}^{k}(\delta ,\lambda ,\mu
,l)}{(B-A)(p-\sigma )}\right) .  \tag{5.2}
\end{equation}%
Then for $m\geq k+p,$ we have%
\begin{equation}
\Re \left( \frac{f(z)}{\kappa _{m}(z)}\right) >1-\frac{1}{{\theta }_{m+1}} 
\tag{5.3}
\end{equation}%
and%
\begin{equation}
\Re \left( \frac{\kappa _{m}(z)}{f(z)}\right) >\frac{{\theta }_{m+1}}{1+{%
\theta }_{m+1}}.  \tag{5.4}
\end{equation}%
The results are sharp for every $m$ with the extremal functions given by%
\begin{equation}
f(z)=z^{p}-\frac{1}{{\theta }_{m+1}}z^{m+1}.  \tag{5.5}
\end{equation}

\noindent \textbf{Proof. }Under the hypothesis of the theorem, we can see
from (5.2) that%
\begin{equation*}
{\theta }_{k+1}>{\theta }_{k}>1\text{\ \ }(k\geq n+p).
\end{equation*}%
Therefore, we have%
\begin{equation}
\sum_{k=n+p}^{m}\left\vert a_{k}\right\vert +{\theta }_{m+1}\sum_{k=m+1}^{%
\infty }\left\vert a_{k}\right\vert \leq \sum_{k=n+p}^{\infty }{\theta }%
_{k}\left\vert a_{k}\right\vert \leq 1  \tag{5.6}
\end{equation}%
by using hypothesis (5.2) again.

Upon setting%
\begin{eqnarray}
{\omega (z)} &=&{{\theta }_{m+1}\left[ \frac{f(z)}{\kappa _{m}(z)}-\left( 1-%
\frac{1}{{\theta }_{m+1}}\right) \right] }  \TCItag{5.7} \\
&{=}&{1-\frac{{\theta }_{m+1}\sum_{k=m+1}^{\infty }\left\vert
a_{k}\right\vert z^{k-p}}{1-\sum_{k=n+p}^{m}\left\vert a_{k}\right\vert
z^{k-p}}.}  \notag
\end{eqnarray}%
By applying (5.6) and (5.7), we find that%
\begin{eqnarray}
{\left\vert \frac{\omega (z)-1}{\omega (z)+1}\right\vert } &=&{\left\vert 
\frac{-{\theta }_{m+1}\sum_{k=m+1}^{\infty }\left\vert a_{k}\right\vert
z^{k-p}}{2-2\sum_{k=n+p}^{m}\left\vert a_{k}\right\vert z^{k-p}-{\theta }%
_{m+1}\sum_{k=m+1}^{\infty }\left\vert a_{k}\right\vert z^{k-p}}\right\vert }
\TCItag{5.8} \\
&\leq &{\frac{{\theta }_{m+1}\sum_{k=m+1}^{\infty }\left\vert
a_{k}\right\vert z^{k-p}}{2-2\sum_{k=n+p}^{m}\left\vert a_{k}\right\vert
z^{k-p}-{\theta }_{m+1}\sum_{k=m+1}^{\infty }\left\vert a_{k}\right\vert
z^{k-p}}\leq 1\;\;(z\in \mathrm{{\mathcal{U}}};\;k\geq n+p),}  \notag
\end{eqnarray}%
which shows that $\Re \left( \omega (z)\right) >0\;(z\in \mathrm{{\mathcal{U}%
}}).$ From (5.7), we immediately obtain the inequality (5.3).

To see that the function $f$ given by (5.5) gives the sharp result, we
observe for $z\rightarrow 1^{-}$ that

\begin{equation*}
\frac{f(z)}{\kappa _{m}(z)}=1-\frac{1}{{\theta }_{m+1}}z^{m-p+1}\rightarrow
1-\frac{1}{{\theta }_{m+1}},
\end{equation*}%
which shows that the bound in (5.3) is the best possible.

Similarly, if we put%
\begin{eqnarray}
{\phi (z)} &=&{(1+{\theta }_{m+1})\left[ \frac{\kappa _{m}(z)}{f(z)}-\frac{{%
\theta }_{m+1}}{1+{\theta }_{m+1}}\right] }  \TCItag{5.9} \\
&{=}&{1+\frac{(1+{\theta }_{m+1})\sum_{k=m+1}^{\infty }\left\vert
a_{k}\right\vert z^{k-p}}{1-\sum_{k=n+p}^{m}\left\vert a_{k}\right\vert
z^{k-p}},}  \notag
\end{eqnarray}%
and make use of (5.6), we can deduce that%
\begin{eqnarray}
{\left\vert \frac{\phi (z)-1}{\phi (z)+1}\right\vert } &=&{\left\vert \frac{%
(1+{{\theta }_{m+1}})\sum_{k=m+1}^{\infty }\left\vert a_{k}\right\vert
z^{k-p}}{2-2\sum_{k=n+p}^{m}\left\vert a_{k}\right\vert z^{k-p}+({{\theta }%
_{m+1}}-1)\sum_{k=m+1}^{\infty }\left\vert a_{k}\right\vert z^{k-p}}%
\right\vert }  \TCItag{5.10} \\
&{\leq }&{\;\frac{(1+{{\theta }_{m+1}})\sum_{k=m+1}^{\infty }\left\vert
a_{k}\right\vert z^{k-p}}{2-2\sum_{k=n+p}^{m}\left\vert a_{k}\right\vert
z^{k-p}-({{\theta }_{m+1}}-1)\sum_{k=m+1}^{\infty }\left\vert
a_{k}\right\vert z^{k-p}}\leq 1\;\;(z\in \mathrm{{\mathcal{U}}};\;k\geq n+p),%
}  \notag
\end{eqnarray}%
which leads us immediately to assertion (5.4) of the theorem.

The bound in (5.4) is sharp with the extremal function given by (5.5). The
proof of theorem is thus completed.

\section{\noindent PROPERTIES ASSOCIATED WITH QUASI-CONVOLUTION\textbf{\ }}

In this part, we establish certain results concerning the Quasi-convolution
of function is in the class $\widetilde{\mathcal{P}}_{\lambda ,\mu
,l}^{\delta }(A,B;\sigma ,p).$

For the functions $f_{j}(z)\in \mathrm{{\mathcal{A}}}(n,p)$ given by%
\begin{equation}
f_{j}(z)=z^{p}-\sum_{k=n+p}^{\infty }\left\vert a_{k,j}\right\vert
z^{k}\;\;(j=\overline{1,m},\;p\in \mathrm{%
\mathbb{N}
}),  \tag{6.1}
\end{equation}%
we denote by $(f_{1}\bullet f_{2})(z)$ the Quasi-convolution of functions $%
f_{1}(z)$ and $f_{2}(z),$ that is,%
\begin{equation}
(f_{1}\bullet f_{2})(z)=z^{p}-\sum_{k=n+p}^{\infty }\left\vert
a_{k,1}\right\vert \left\vert a_{k,2}\right\vert z^{k}.  \tag{6.2}
\end{equation}

\noindent \textbf{Theorem 6.1. }If $f_{j}(z)\in \widetilde{\mathcal{P}}%
_{\lambda ,\mu ,l}^{\delta }(A,B;\sigma _{j},p)$ $(j=\overline{1,m}),$ then%
\begin{equation}
(f_{1}\bullet f_{2}\bullet ...\bullet f_{m})(z)\in \widetilde{\mathcal{P}}%
_{\lambda ,\mu ,l}^{\delta }(A,B;\Upsilon ,p),  \tag{6.3}
\end{equation}%
where%
\begin{equation}
\Upsilon :=p-\frac{\prod_{j=1}^{m}(B-A)(p-\sigma _{j})}{(B-A)[(n+p)(1+B)\Phi
_{p}^{n+p}(\delta ,\lambda ,\mu ,l)]^{m-1}}.  \tag{6.4}
\end{equation}%
The result is sharp for the functions $f_{j}(z)$ given by%
\begin{equation}
f_{j}(z)=z^{p}-\frac{(B-A)(p-\sigma _{j})}{(n+p)(1+B)\Phi _{p}^{n+p}(\delta
,\lambda ,\mu ,l)}z^{p+n}\;\;(j=\overline{1,m}).  \tag{6.5}
\end{equation}

\noindent

\noindent \textbf{Proof. }For $m=1,$ we see that $\Upsilon =\sigma _{1}.$
For $m=2,$ Theorem 3.1 gives%
\begin{equation}
\sum_{k=n+p}^{\infty }\frac{k(1+B)\Phi _{p}^{k}(\delta ,\lambda ,\mu ,l)}{%
(B-A)(p-\sigma _{j})}\left\vert a_{k,j}\right\vert \leq 1\;\;(j=1,2). 
\tag{6.6}
\end{equation}%
Therefore, by the Cauchy-Schwarz inequality, we obtain%
\begin{equation}
\sum_{k=n+p}^{\infty }\frac{k(1+B)\Phi _{p}^{k}(\delta ,\lambda ,\mu ,l)}{%
\sqrt{\prod_{j=1}^{2}(B-A)(p-\sigma _{j})}}\sqrt{\left\vert
a_{k,1}\right\vert \left\vert a_{k,2}\right\vert }\leq 1.  \tag{6.7}
\end{equation}%
To prove the case when $m=2,$ we have to find the largest $\Upsilon $ such
that%
\begin{equation}
\sum_{k=n+p}^{\infty }\frac{k(1+B)\Phi _{p}^{k}(\delta ,\lambda ,\mu ,l)}{%
(B-A)(p-\Upsilon )}\left\vert a_{k,1}\right\vert \left\vert
a_{k,2}\right\vert \leq 1,  \tag{6.8}
\end{equation}%
or such that%
\begin{equation}
\frac{\left\vert a_{k,1}\right\vert \left\vert a_{k,2}\right\vert }{%
(B-A)(p-\Upsilon )}\leq \frac{\sqrt{\left\vert a_{k,1}\right\vert \left\vert
a_{k,2}\right\vert }}{\sqrt{\prod_{j=1}^{2}(B-A)(p-\sigma _{j})}},  \tag{6.9}
\end{equation}%
this, equivalently, that%
\begin{equation}
\sqrt{\left\vert a_{k,1}\right\vert \left\vert a_{k,2}\right\vert }\leq 
\frac{(B-A)(p-\Upsilon )}{\sqrt{\prod_{j=1}^{2}(B-A)(p-\sigma _{j})}}. 
\tag{6.10}
\end{equation}%
Further, by using (6.7), we need to find the largest $\Upsilon $ such that%
\begin{equation*}
\frac{\sqrt{\prod_{j=1}^{2}(B-A)(p-\sigma _{j})}}{k(1+B)\Phi _{p}^{k}(\delta
,\lambda ,\mu ,l)}\leq \frac{(B-A)(p-\Upsilon )}{\sqrt{%
\prod_{j=1}^{2}(B-A)(p-\sigma _{j})}}
\end{equation*}%
or, equivalently, that%
\begin{equation}
\frac{1}{(B-A)(p-\Upsilon )}\leq \frac{k(1+B)\Phi _{p}^{k}(\delta ,\lambda
,\mu ,l)}{\prod_{j=1}^{2}(B-A)(p-\sigma _{j})}.  \tag{6.11}
\end{equation}%
It follows from (6.9) that%
\begin{equation}
\Upsilon \leq p-\frac{\prod_{j=1}^{2}(B-A)(p-\sigma _{j})}{(B-A)k(1+B)\Phi
_{p}^{k}(\delta ,\lambda ,\mu ,l)}.  \tag{6.12}
\end{equation}%
Now, defining the function $\psi (k)$ by%
\begin{equation}
\psi (k)=p-\frac{\prod_{j=1}^{2}(B-A)(p-\sigma _{j})}{(B-A)k(1+B)\Phi
_{p}^{k}(\delta ,\lambda ,\mu ,l)},  \tag{6.13}
\end{equation}%
we see that $\psi ^{\prime }(k)\geq 0$ for $k\geq p+n.$ This implies that%
\begin{equation*}
\Upsilon \leq \psi (n+p)=p-\frac{\prod_{j=1}^{2}(B-A)(p-\sigma _{j})}{%
(B-A)(n+p)(1+B)\Phi _{p}^{n+p}(\delta ,\lambda ,\mu ,l)}.
\end{equation*}%
Therefore, the result is true for $m=2.$

Suppose that the result is true for any positive integer $m.$ Then we have $%
(f_{1}\bullet f_{2}\bullet ...\bullet f_{m}\bullet f_{m+1})(z)\in \widetilde{%
\mathcal{P}}_{\lambda ,\mu ,l}^{\delta }(A,B;\gamma ,p),$ when%
\begin{equation*}
\gamma =p-\frac{(B-A)(p-\Upsilon )(B-A)(p-\sigma _{m+1})}{%
(B-A)(n+p)(1+B)\Phi _{p}^{n+p}(\delta ,\lambda ,\mu ,l)}
\end{equation*}%
where $\Upsilon $ is given by (6.4). After a simple calculation, we have%
\begin{equation*}
\gamma \leq p-\frac{\prod_{j=1}^{m+1}(B-A)(p-\sigma _{j})}{%
(B-A)[(n+p)(1+B)\Phi _{p}^{n+p}(\delta ,\lambda ,\mu ,l)]^{m}}.
\end{equation*}%
Thus, the result is true for $m+1.$ Therefore, by using the mathematical
induction, we conclude that the result is true for any positive integer $m.$

Finally, taking the functions $f_{j}(z)$ defined by (6.5), we have%
\begin{eqnarray*}
{(f_{1}\bullet f_{2}\bullet ...\bullet f_{m})(z)} &=&{z^{p}-\left\{
\prod_{j=1}^{m}\frac{(B-A)(p-\sigma _{j})}{(p+n)(1+B)\Phi _{p}^{p+n}(\delta
,\lambda ,\mu ,l)}\right\} z^{p+n}} \\
&{=}&{z^{p}-\mathrm{A}_{p+n}z^{p+n},}
\end{eqnarray*}%
which shows that%
\begin{eqnarray*}
~~~~~\sum_{k=p+n}^{\infty }\frac{k(1+B)\Phi _{p}^{k}(\delta ,\lambda ,\mu ,l)%
}{(B-A)(p-\Upsilon )}\mathrm{A}_{k} &=&\frac{(n+p)(1+B)\Phi
_{p}^{n+p}(\delta ,\lambda ,\mu ,l)}{(B-A)(p-\Upsilon )}\mathrm{A}%
_{p+n}~~~~~~~~~~~~~ \\
~~~~~~~~~~~~~ &=&\frac{(n+p)(1+B)\Phi _{p}^{n+p}(\delta ,\lambda ,\mu ,l)}{%
(B-A)(p-\Upsilon )} \\
&&\times \left\{ \prod_{j=1}^{2}\frac{(B-A)(p-\sigma _{j})}{(p+n)(1+B)\Phi
_{p}^{p+n}(\delta ,\lambda ,\mu ,l)}\right\} .
\end{eqnarray*}%
Consequently, the result is sharp.

Putting $\sigma _{j}=\sigma $ $(j=\overline{1,m})$ in Theorem 6.1, we have;

\noindent \textbf{Corollary 6.2. }If $f_{j}(z)\in \widetilde{\mathcal{P}}%
_{\lambda ,\mu ,l}^{\delta }(A,B;\sigma ,p)$ $(j=\overline{1,m}),$ then%
\begin{equation*}
(f_{1}\bullet f_{2}\bullet ...\bullet f_{m})(z)\in \widetilde{\mathcal{P}}%
_{\lambda ,\mu ,l}^{\delta }(A,B;\Upsilon ,p),
\end{equation*}%
where%
\begin{equation*}
\Upsilon :=p-\frac{[(B-A)(p-\sigma )]^{m}}{(B-A)[(n+p)(1+B)\Phi
_{p}^{n+p}(\delta ,\lambda ,\mu ,l)]^{m-1}}.
\end{equation*}%
The result is sharp for the functions $f_{j}(z)$ given by%
\begin{equation*}
f_{j}(z)=\frac{(B-A)(p-\sigma )}{(n+p)(1+B)\Phi _{p}^{n+p}(\delta ,\lambda
,\mu ,l)}z^{p+n}\;\;(j=\overline{1,m}).
\end{equation*}

\noindent

\textbf{Remark 6.1. }For special values of parameters $\lambda ,\mu
,l,\delta ,\sigma ,A,B,n$ and $p,$ our results reduce to several well-known
results as follows:

\begin{enumerate}
\item Putting $A=-1,\;B=1,\;\delta =0$ and $m=2$ in Theorem 6.1, we obtain
the corresponding results of Yaguchi \textit{et al}. [55] and Aouf and
Darwish [7] for $n=1.$

\item Putting $A=-1,\;B=1,\;\delta =0$ and $m=2$ in Corollary 6.2, we obtain
the corresponding results of Lee \textit{et al}. [26] and for $n=1$ and
Sekine and Owa [44] for $p=1.$

\item Putting $A=-1,\;B=1,\;\delta =0$ and $m=3$ in Corollary 6.2, we obtain
the corresponding result due to Aouf and Darwish [7] for $n=1.$

\item Putting $A=-\beta ,\;B=\beta ,\;\delta =0$ and $\sigma =\alpha $ in
Theorem 6.1, we obtain the corresponding result due to Aouf [5].
\end{enumerate}

\noindent \textbf{Theorem 6.2. }Let the function $f_{j}(z)$ $(j=\overline{1,m%
})$ given by (6.1) be in the class $\widetilde{\mathcal{P}}_{\lambda ,\mu
,l}^{\delta }(A,B;\sigma _{j},p).$ Then the function%
\begin{equation}
h(z)=z^{p}-\sum_{k=n+p}^{\infty }\left( \sum_{j=1}^{m}\left\vert
a_{k,j}\right\vert ^{2}\right) z^{k}  \tag{6.14}
\end{equation}%
belongs to the class $\widetilde{\mathcal{P}}_{\lambda ,\mu ,l}^{\delta
}(A,B;\chi ,p),$ where%
\begin{equation}
\chi :=p-\frac{m(B-A)(p-\sigma ^{\ast })^{2}}{(n+p)(1+B)\Phi
_{p}^{n+p}(\delta ,\lambda ,\mu ,l)}\;\;(\sigma ^{\ast }:=\min \{\sigma
_{1},\sigma _{2},...,\sigma _{m}\}).  \tag{6.15}
\end{equation}%
The result is sharp for the functions $f_{j}(z)$ $(j=\overline{1,m})$ given
by (6.5).

\noindent \textbf{Proof. }By virtue of Theorem 3.1 we have%
\begin{equation}
\sum_{k=p+n}^{\infty }\left\{ \frac{k(1+B)\Phi _{p}^{k}(\delta ,\lambda ,\mu
,l)}{(B-A)(p-\sigma _{j})}\right\} ^{2}\left\vert a_{k,j}\right\vert
^{2}\leq \left\{ \sum_{k=p+n}^{\infty }\frac{k(1+B)\Phi _{p}^{k}(\delta
,\lambda ,\mu ,l)}{(B-A)(p-\sigma _{j})}\left\vert a_{k,j}\right\vert
\right\} ^{2}\leq 1.  \tag{6.16}
\end{equation}%
Then it follows that for $(j=\overline{1,m}),$%
\begin{equation}
\frac{1}{m}\sum_{k=p+n}^{\infty }\left\{ \frac{k(1+B)\Phi _{p}^{k}(\delta
,\lambda ,\mu ,l)}{(B-A)(p-\sigma _{j})}\right\} ^{2}\left(
\sum_{j=1}^{m}\left\vert a_{k,j}\right\vert ^{2}\right) \leq 1.  \tag{6.17}
\end{equation}%
Therefore, we need to find the largest $\chi $ such that%
\begin{equation}
\frac{1}{m}\sum_{k=p+n}^{\infty }\left\{ \frac{k(1+B)\Phi _{p}^{k}(\delta
,\lambda ,\mu ,l)}{(B-A)(p-\chi )}\right\} \left( \sum_{j=1}^{m}\left\vert
a_{k,j}\right\vert ^{2}\right) \leq 1.  \tag{6.18}
\end{equation}%
This implies that%
\begin{equation}
\chi \leq p-\frac{m(B-A)(p-\sigma ^{\ast })^{2}}{k(1+B)\Phi _{p}^{k}(\delta
,\lambda ,\mu ,l)}\;\;(\sigma ^{\ast }:=\min \{\sigma _{1},\sigma
_{2},...,\sigma _{m}\},\;k\geq p+n).  \tag{6.19}
\end{equation}%
Now, defining the function $\Im (k)$ by%
\begin{equation}
\Im (k):=p-\frac{m(B-A)(p-\sigma ^{\ast })^{2}}{k(1+B)\Phi _{p}^{k}(\delta
,\lambda ,\mu ,l)},  \tag{6.20}
\end{equation}%
we see that $\Im (k)$ is an increasing function of $k,$ $k\geq p+n.$ Setting 
$k=p+n$ in (6.19) we have%
\begin{equation*}
\chi \leq \Im (n+p):=p-\frac{m(B-A)(p-\sigma ^{\ast })^{2}}{(n+p)(1+B)\Phi
_{p}^{n+p}(\delta ,\lambda ,\mu ,l)}
\end{equation*}%
which completes the proof of Theorem 6.2.

Setting $\sigma _{j}=\sigma $ $(j=\overline{1,m}),$ in Theorem 6.2, we
arrive at the following result.

\noindent \textbf{Corollary 6.3. }Let the functions $f_{j}(z)$ $(j=\overline{%
1,m})$ given by (6.1) be in the class $\widetilde{\mathcal{P}}_{\lambda ,\mu
,l}^{\delta }(A,B;\sigma ,p).$ Then the function%
\begin{equation*}
h(z)=z^{p}-\sum_{k=n+p}^{\infty }\left( \sum_{j=1}^{m}\left\vert
a_{k,j}\right\vert ^{2}\right) z^{k}
\end{equation*}%
belongs to the class $\widetilde{\mathcal{P}}_{\lambda ,\mu ,l}^{\delta
}(A,B;\chi ,p),$ where%
\begin{equation*}
\chi :=p-\frac{m(B-A)(p-\sigma )^{2}}{(n+p)(1+B)\Phi _{p}^{n+p}(\delta
,\lambda ,\mu ,l)}.
\end{equation*}%
The result is sharp for the functions $f_{j}(z)$ $(j=\overline{1,m})$ given
by (6.5).

\noindent \textbf{Remark 6.2.}

\begin{enumerate}
\item Putting $A=-1,\;B=1,\;\delta =0$ and $m=2$ in Theorem 6.2, we obtain
the corresponding results of Yaguchi \textit{et al}. [55].

\item Putting $A=-1,\;B=1,\;\delta =0$ and $m=2$ in Corollary 6.3, we obtain
the corresponding results of Aouf and Darwish [7] for $n=1,$ Sekine and Owa
[44] for $p=1$.

\item Putting $A=-\beta ,\;B=\beta ,\;\delta =0$ and $\sigma =\alpha $ in
Theorem 6.2, we obtain the corresponding result due to Aouf [5].
\end{enumerate}

\section{\noindent APPLICATIONS OF FRACTIONAL CALCULUS OPERATORS}

Various operators of fractional calculus (that is, fractional integral and
fractional derivatives) have been studied in the literature rather
extensively (\textit{cf., e.g}., [37, 50, 53]; see also [14, 51] the various
references cited therein). For our present investigation, we recall the
following definitions.

\noindent \textbf{Definition 7.1.} Let $f(z)$ be analytic in a simply
connected region of the $z$-\textit{plane} containing the origin. The
fractional integral of $f$ of order $\nu $ is defined by%
\begin{equation}
\mathrm{{\mathcal{D}}}_{z}^{-\nu }f(z)=\frac{1}{\Gamma (\nu )}\int_{0}^{z}%
\frac{f(\zeta )}{(z-\zeta )^{1-\nu }}d\zeta \text{\ \ }(\nu >0),  \tag{7.1}
\end{equation}%
where the multiplicity of $(z-\zeta )^{\nu -1}$ is removed by requiring that 
$\log (z-\zeta )$ is real for $z-\zeta >0$.

\noindent

\noindent \textbf{Definition 7.2.} Let $f(z)$ be analytic in a simply
connected region of the $z$-\textit{plane} containing the origin. The
fractional derivative of $f$ of order $\nu $ is defined by%
\begin{equation}
\mathrm{{\mathcal{D}}}_{z}^{\nu }f(z)=\frac{1}{\Gamma (1-\nu )}\int_{0}^{z}%
\frac{f(\zeta )}{(z-\zeta )^{\nu }}d\zeta \text{\ \ }(0\leq \nu <1), 
\tag{7.2}
\end{equation}%
where the multiplicity of $(z-\zeta )^{-\nu }$ is removed by requiring that $%
\log (z-\zeta )$ is real for $z-\zeta >0$.

\noindent

\noindent \textbf{Definition 7.3}. Under the hypotheses of Definition 7.2,
the fractional derivative of order $n+\nu $ is defined, for a function $f(z)$%
, by%
\begin{equation}
\mathrm{{\mathcal{D}}}_{z}^{n+\nu }f(z)=\frac{d^{n}}{dz^{n}}\{\mathrm{{%
\mathcal{D}}}_{z}^{\nu }f(z)\}\;\;(0\leq \nu <1;\;n\in \mathrm{%
\mathbb{N}
}_{0}).  \tag{7.3}
\end{equation}%
In this section, we shall investigate the growth and distortion properties
of functions in the class $\widetilde{\mathcal{P}}_{\lambda ,\mu ,l}^{\delta
}(A,B;\sigma ,p),$ which involving the operators $\mathrm{{\mathcal{I}}}%
_{\vartheta ,p}$ and $\mathrm{{\mathcal{D}}}_{z}^{\nu }.$ In order to derive
our results, we need the following lemma given by Chen \textit{et al.} [14].

\noindent \textbf{Lemma 7.1 }(\textit{see} [14]). Let the function $f(z)$
defined by (1.10). Then%
\begin{equation}
\mathrm{{\mathcal{D}}}_{z}^{\nu }\{(\mathrm{{\mathcal{I}}}_{\vartheta
,p}f)(z)\}=\frac{\Gamma (p+1)}{\Gamma (p+1-\nu )}z^{p-\nu
}-\sum_{k=n+p}^{\infty }\frac{(\vartheta +p)\Gamma (k+1)}{(\vartheta
+k)\Gamma (k+1-\nu )}a_{k}z^{k-\nu }  \tag{7.4}
\end{equation}%
\begin{equation*}
(\nu \in \mathrm{%
\mathbb{R}
};\;\vartheta >-p;\;p,n\in \mathrm{%
\mathbb{N}
})
\end{equation*}%
and%
\begin{equation}
\mathrm{{\mathcal{I}}}_{\vartheta ,p}\{(\mathrm{{\mathcal{D}}}_{z}^{\nu
}f)(z)\}=\frac{(\vartheta +p)\Gamma (p+1)}{(\vartheta +p-\nu )\Gamma
(p+1-\nu )}z^{p-\nu }-\sum_{k=n+p}^{\infty }\frac{(\vartheta +p)\Gamma (k+1)%
}{(\vartheta +k-\nu )\Gamma (k+1-\nu )}a_{k}z^{k-\nu }  \tag{7.5}
\end{equation}%
\begin{equation*}
(\nu \in \mathrm{%
\mathbb{R}
};\;\vartheta >-p;\;p,n\in \mathrm{%
\mathbb{N}
})
\end{equation*}%
provided that no zeros appear in the denominators in (7.4) and (7.5).

\noindent \textbf{Theorem 7.1. }Let the functions $f(z)$ defined by (1.10)
be in the class $\widetilde{\mathcal{P}}_{\lambda ,\mu ,l}^{\delta
}(A,B;\sigma ,p).$ Then%
\begin{eqnarray}
&&{\left\vert \mathrm{{\mathcal{D}}}_{z}^{-\nu }\{(\mathrm{{\mathcal{I}}}%
_{\vartheta ,p}f)(z)\}\right\vert }  \TCItag{7.6} \\
&{\geq }&\left\{ {\frac{\Gamma (p+1)}{\Gamma (p+1+\nu )}-\frac{(\vartheta
+p)\Gamma (n+p+1)(B-A)(p-\sigma )}{(\vartheta +n+p)\Gamma (n+p+1+\nu
)(n+p)(1+B)\Phi _{p}^{n+p}(\delta ,\lambda ,\mu ,l)}\left\vert z\right\vert
^{n}}\right\} {\left\vert z\right\vert ^{p+\nu }}  \notag
\end{eqnarray}%
\begin{equation*}
(z\in \mathrm{{\mathcal{U}}};\;\nu >0;\;\vartheta >-p;\;p,n\in \mathrm{%
\mathbb{N}
})
\end{equation*}%
and%
\begin{eqnarray}
&&{\left\vert \mathrm{{\mathcal{D}}}_{z}^{-\nu }\{(\mathrm{{\mathcal{I}}}%
_{\vartheta ,p}f)(z)\}\right\vert }  \TCItag{7.7} \\
&{\leq }&\left\{ {\frac{\Gamma (p+1)}{\Gamma (p+1+\nu )}+\frac{(\vartheta
+p)\Gamma (n+p+1)(B-A)(p-\sigma )}{(\vartheta +n+p)\Gamma (n+p+1+\nu
)(n+p)(1+B)\Phi _{p}^{n+p}(\delta ,\lambda ,\mu ,l)}\left\vert z\right\vert
^{n}}\right\} {\left\vert z\right\vert ^{p+\nu }}  \notag
\end{eqnarray}%
\begin{equation*}
(z\in \mathrm{{\mathcal{U}}};\;\nu >0;\;\vartheta >-p;\;p,n\in \mathrm{%
\mathbb{N}
}).
\end{equation*}%
Each of the assertions (7.6) and (7.7) is sharp.

\noindent \textbf{Proof.} In view of Theorem 3.1, we have%
\begin{equation}
\frac{(n+p)(1+B)\Phi _{p}^{n+p}(\delta ,\lambda ,\mu ,l)}{(B-A)(p-\sigma )}%
\sum_{k=n+p}^{\infty }\left\vert a_{k}\right\vert \leq \sum_{k=n+p}^{\infty }%
\frac{k(1+B)\Phi _{p}^{k}(\delta ,\lambda ,\mu ,l)}{(B-A)(p-\sigma )}%
\left\vert a_{k}\right\vert \leq 1,  \tag{7.8}
\end{equation}%
which readily yields%
\begin{equation}
\sum_{k=n+p}^{\infty }\left\vert a_{k}\right\vert \leq \frac{(B-A)(p-\sigma )%
}{(n+p)(1+B)\Phi _{p}^{n+p}(\delta ,\lambda ,\mu ,l)}.  \tag{7.9}
\end{equation}%
Consider the function $\mathrm{{\mathcal{F}}}(z)$ defined in $\mathrm{{%
\mathcal{U}}}$ by%
\begin{eqnarray*}
\mathrm{{\mathcal{F}}}(z) &=&\frac{\Gamma (p+1-\nu )}{\Gamma (p+1)}z^{-\nu }%
\mathrm{{\mathcal{D}}}_{z}^{\nu }\{(\mathrm{{\mathcal{I}}}_{\vartheta
,p}f)(z)\} \\
&=&z^{p}-\sum_{k=n+p}^{\infty }\frac{(\vartheta +p)\Gamma (k+1)\Gamma
(p+1+\nu )}{(\vartheta +k)\Gamma (k+1+\nu )\Gamma (p+1)}\left\vert
a_{k}\right\vert z^{k} \\
&=&z^{p}-\sum_{k=n+p}^{\infty }\Theta (k)\left\vert a_{k}\right\vert z^{k}%
\text{ \ }(z\in \mathrm{{\mathcal{U}}})
\end{eqnarray*}%
where

\begin{equation}
\Theta (k):=\frac{(\vartheta +p)\Gamma (k+1)\Gamma (p+1+\nu )}{(\vartheta
+k)\Gamma (k+1+\nu )\Gamma (p+1)}\text{\ \ }(k\geq p+n;\;p,n\in \mathrm{%
\mathbb{N}
};\;\nu >0).  \tag{7.10}
\end{equation}%
Since $\Theta (k)$ is a \textit{decreasing} function of $k$ when $\nu >0,$
we get%
\begin{equation}
0<\Theta (k)\leq \Theta (n+p)=\frac{(\vartheta +p)\Gamma (n+p+1)\Gamma
(p+1+\nu )}{(\vartheta +n+p)\Gamma (n+p+1+\nu )\Gamma (p+1)}  \tag{7.11}
\end{equation}%
\begin{equation*}
(\nu >0;\;\vartheta >-p;\;p,n\in \mathrm{%
\mathbb{N}
}).
\end{equation*}%
Thus, by using (7.9) and (7.11), for all $z\in \mathcal{U},$ we deduce that%
\begin{eqnarray*}
\left\vert \mathrm{{\mathcal{F}}}(z)\right\vert &\geq &\left\vert
z\right\vert ^{p}-\Theta (n+p)\left\vert z\right\vert
^{n+p}\sum_{k=n+p}^{\infty }\left\vert a_{k}\right\vert \\
&\geq &\left\vert z\right\vert ^{p}-\frac{(\vartheta +p)\Gamma (n+p+1)\Gamma
(p+1+\nu )(B-A)(p-\sigma )}{(\vartheta +n+p)\Gamma (n+p+1+\nu )\Gamma
(p+1)(n+p)(1+B)\Phi _{p}^{n+p}(\delta ,\lambda ,\mu ,l)}\left\vert
z\right\vert ^{n+p}
\end{eqnarray*}%
and%
\begin{eqnarray*}
\left\vert \mathrm{{\mathcal{F}}}(z)\right\vert &\leq &\left\vert
z\right\vert ^{p}+\Theta (n+p)\left\vert z\right\vert
^{n+p}\sum_{k=n+p}^{\infty }\left\vert a_{k}\right\vert \\
&\leq &\left\vert z\right\vert ^{p}+\frac{(\vartheta +p)\Gamma (n+p+1)\Gamma
(p+1+\nu )(B-A)(p-\sigma )}{(\vartheta +n+p)\Gamma (n+p+1+\nu )\Gamma
(p+1)(n+p)(1+B)\Phi _{p}^{n+p}(\delta ,\lambda ,\mu ,l)}\left\vert
z\right\vert ^{n+p}
\end{eqnarray*}%
which yield the inequalities (7.6) and (7.7) of Theorem 7.1. Equalities in
(7.6) and (7.7) are attained for the function $f(z)$ given by

\begin{eqnarray*}
{\mathrm{{\mathcal{D}}}_{z}^{-\nu }\{(\mathrm{{\mathcal{I}}}_{\vartheta
,p}f)(z)\}} &=&{\left\{ \frac{\Gamma (p+1)}{\Gamma (p+1+\nu )}\right. } \\
&&{\,\,\left. \,-\frac{(\vartheta +p)\Gamma (n+p+1)(B-A)(p-\sigma )}{%
(\vartheta +n+p)\Gamma (n+p+1+\nu )(n+p)(1+B)\Phi _{p}^{n+p}(\delta ,\lambda
,\mu ,l)}z^{n}\right\} z^{p+\nu }}
\end{eqnarray*}%
or, equivalently, by%
\begin{equation*}
(\mathrm{{\mathcal{I}}}_{\vartheta ,p}f)(z)=z^{p}-\frac{(\vartheta
+p)(B-A)(p-\sigma )}{(\vartheta +n+p)(n+p)(1+B)\Phi _{p}^{n+p}(\delta
,\lambda ,\mu ,l)}z^{n+p}.
\end{equation*}%
Thus, we complete the proof of Theorem 7.1.

\noindent \textbf{Theorem 7.2. }Let the functions $f(z)$ defined by (1.10)
be in the class $\widetilde{\mathcal{P}}_{\lambda ,\mu ,l}^{\delta
}(A,B;\sigma ,p).$ Then%
\begin{eqnarray}
&&{\left\vert \mathrm{{\mathcal{D}}}_{z}^{\nu }\{(\mathrm{{\mathcal{I}}}%
_{\vartheta ,p}f)(z)\}\right\vert }  \TCItag{7.12} \\
&{\geq }&\left\{ {\frac{\Gamma (p+1)}{\Gamma (p+1-\nu )}z^{p-\nu }-\frac{%
(\vartheta +p)\Gamma (n+p)(B-A)(p-\sigma )}{(\vartheta +n+p)\Gamma
(n+p+1-\nu )(1+B)\Phi _{p}^{n+p}(\delta ,\lambda ,\mu ,l)}\left\vert
z\right\vert ^{n}}\right\} {\left\vert z\right\vert ^{p-\nu }}  \notag
\end{eqnarray}%
\begin{equation*}
(z\in \mathrm{{\mathcal{U}}};\;\nu >0;\;\vartheta >-p;\;p,n\in \mathrm{%
\mathbb{N}
})
\end{equation*}%
and%
\begin{eqnarray}
&&{\left\vert \mathrm{{\mathcal{D}}}_{z}^{\nu }\{(\mathrm{{\mathcal{I}}}%
_{\vartheta ,p}f)(z)\}\right\vert }  \TCItag{7.13} \\
&{\leq }&\left\{ {\frac{\Gamma (p+1)}{\Gamma (p+1-\nu )}z^{p-\nu }+\frac{%
(\vartheta +p)\Gamma (n+p)(B-A)(p-\sigma )}{(\vartheta +n+p)\Gamma
(n+p+1-\nu )(1+B)\Phi _{p}^{n+p}(\delta ,\lambda ,\mu ,l)}\left\vert
z\right\vert ^{n}}\right\} {\left\vert z\right\vert ^{p-\nu }}  \notag
\end{eqnarray}

\begin{equation*}
(z\in \mathrm{{\mathcal{U}}};\;\nu >0;\;\vartheta >-p;\;p,n\in \mathrm{%
\mathbb{N}
}).
\end{equation*}%
Each of the assertions (7.12) and (7.13) is sharp.

\noindent \textbf{Proof.} It follows from Theorem 3.1 that\textbf{\ }%
\begin{equation}
\sum_{k=n+p}^{\infty }k\left\vert a_{k}\right\vert \leq \frac{(B-A)(p-\sigma
)}{(1+B)\Phi _{p}^{n+p}(\delta ,\lambda ,\mu ,l)}.  \tag{7.14}
\end{equation}%
Consider the function $\mathrm{{\mathcal{Q}}}(z)$ defined in $\mathrm{{%
\mathcal{U}}}$ by%
\begin{eqnarray*}
\mathrm{{\mathcal{Q}}}(z) &=&\frac{\Gamma (p+1-\nu )}{\Gamma (p+1)}z^{\nu }%
\mathrm{{\mathcal{D}}}_{z}^{\nu }\{(\mathrm{{\mathcal{I}}}_{\vartheta
,p}f)(z)\} \\
&=&z^{p}-\sum_{k=n+p}^{\infty }\frac{(\vartheta +p)\Gamma (k)\Gamma (p+1-\nu
)}{(\vartheta +k)\Gamma (k+1-\nu )\Gamma (p+1)}k\left\vert a_{k}\right\vert
z^{k} \\
&=&z^{p}-\sum_{k=n+p}^{\infty }\wp (k)k\left\vert a_{k}\right\vert z^{k}%
\text{\ \ }(z\in \mathrm{{\mathcal{U}}})
\end{eqnarray*}%
where, for convenience,

\begin{equation}
\wp (k):=\frac{(\vartheta +p)\Gamma (k)\Gamma (p+1-\nu )}{(\vartheta
+k)\Gamma (k+1-\nu )\Gamma (p+1)}\text{\ \ }(k\geq p+n;\;p,n\in 
\mathbb{N}
;\;0\leq \nu <1).  \tag{7.15}
\end{equation}%
Since $\wp (k)$ is a \textit{decreasing} function of $k$ when $0\leq \nu <1,$
we find that%
\begin{equation}
0<\wp (k)\leq \wp (n+p)=\frac{(\vartheta +p)\Gamma (n+p)\Gamma (p+1-\nu )}{%
(\vartheta +n+p)\Gamma (n+p+1-\nu )\Gamma (p+1)}  \tag{7.16}
\end{equation}%
\begin{equation*}
(0\leq \nu <1;\;\vartheta >-p;\;p,n\in \mathrm{%
\mathbb{N}
}).
\end{equation*}%
Hence, with the aid of (7.14) and (7.16), for all $z\in \mathcal{U}$\textbf{$%
,$} we have%
\begin{eqnarray*}
\left\vert \mathrm{{\mathcal{Q}}}(z)\right\vert &\geq &\left\vert
z\right\vert ^{p}-\wp (n+p)\left\vert z\right\vert
^{n+p}\sum_{k=n+p}^{\infty }k\left\vert a_{k}\right\vert \\
&\geq &\left\vert z\right\vert ^{p}-\frac{(\vartheta +p)\Gamma (n+p)\Gamma
(p+1-\nu )(B-A)(p-\sigma )}{(\vartheta +n+p)\Gamma (n+p+1-\nu )\Gamma
(p+1)(1+B)\Phi _{p}^{n+p}(\delta ,\lambda ,\mu ,l)}\left\vert z\right\vert
^{n+p}
\end{eqnarray*}%
and%
\begin{eqnarray*}
\left\vert \mathrm{{\mathcal{Q}}}(z)\right\vert &\leq &\left\vert
z\right\vert ^{p}+\wp (n+p)\left\vert z\right\vert
^{n+p}\sum_{k=n+p}^{\infty }k\left\vert a_{k}\right\vert \\
&\leq &\left\vert z\right\vert ^{p}+\frac{(\vartheta +p)\Gamma (n+p)\Gamma
(p+1-\nu )(B-A)(p-\sigma )}{(\vartheta +n+p)\Gamma (n+p+1-\nu )\Gamma
(p+1)(1+B)\Phi _{p}^{n+p}(\delta ,\lambda ,\mu ,l)}\left\vert z\right\vert
^{n+p}
\end{eqnarray*}%
which yield the inequalities (7.15) and (7.16) of Theorem 7.2. Equalities in
(7.15) and (7.16) are attained for the function $f(z)$ given by%
\begin{eqnarray*}
&&{\mathrm{{\mathcal{D}}}_{z}^{\nu }\{(\mathrm{{\mathcal{I}}}_{\vartheta
,p}f)(z)\}} \\
&{=}&\left\{ {\frac{\Gamma (p+1)}{\Gamma (p+1-\nu )}-\frac{(\vartheta
+p)\Gamma (n+p+1)(B-A)(p-\sigma )}{(\vartheta +n+p)\Gamma (n+p+1-\nu
)(1+B)\Phi _{p}^{n+p}(\delta ,\lambda ,\mu ,l)}z^{n}}\right\} {z^{p-\nu }}
\end{eqnarray*}%
or, equivalently, by%
\begin{equation*}
(\mathrm{{\mathcal{I}}}_{\vartheta ,p}f)(z)=z^{p}-\frac{(\vartheta
+p)(B-A)(p-\sigma )}{(\vartheta +n+p)(n+p)(1+B)\Phi _{p}^{n+p}(\delta
,\lambda ,\mu ,l)}z^{n+p}.
\end{equation*}%
Consequently, we complete the proof of Theorem 7.2.

\noindent

\noindent

\end{document}